\documentclass[a4paper,10pt]{scrartcl}
\usepackage[english]{babel}
\usepackage[utf8]{inputenc}
\usepackage{enumerate}
\usepackage{geometry}
\usepackage{amsfonts,amsmath,amssymb,amsopn}
\usepackage{amsthm}
\usepackage{mathtools}
\usepackage{stmaryrd}
\usepackage{nicefrac}
\usepackage{exscale}
\usepackage{relsize}

\usepackage[numbers,square]{natbib}
\usepackage{url} 
\usepackage[colorlinks=true,pdfpagelabels,unicode]{hyperref}
\hypersetup{linkcolor=blue, urlcolor=blue, citecolor=blue} 

\allowdisplaybreaks 

\theoremstyle{plain}

\newtheoremstyle{theo}
	{3pt} 
	{3pt} 
	{\itshape} 
	{} 
		{\bfseries} 
	{\\} 
	{ } 
	{\thmname{#1}\thmnumber{ #2.}\thmnote{ - #3}} 
\theoremstyle{theo}

\newtheorem{definition}{Definition}[section]
\newtheorem{lemma}[definition]{Lemma}
\newtheorem{theorem}[definition]{Theorem}

\newenvironment{bew}{\begin{proof}[\bfseries Proof:]}{\end{proof}}

\DeclareMathOperator{\bomega}{\overline{\Omega}}
\DeclareMathOperator{\romega}{\partial\Omega}
\DeclareMathOperator{\intd}{d\!}
\DeclareMathOperator{\wto}{\rightharpoonup}
\DeclareMathOperator{\wsto}{\stackrel{\star}{\wto}}
\newcommand{\epsi}{\varepsilon}

\usepackage[usenames,dvipsnames,svgnames,table]{xcolor}

\newcommand{\Tmek}{T^{(\kappa)}_{max,\;\!\epsi}}
\newcommand{\GNI}{Gagliardo--Nirenberg inequality}
\newcommand{\CSI}{Cauchy--Schwarz inequality}
\newcommand{\intttau}{\int_t^{t+\tau}\!} 
\newcommand{\into}[1]{\int_0^{#1}\!}

\newcommand{\intoT}{\into{T}}

\newcommand{\intomega}{\int_{\Omega}\!} 
\newcommand{\inttauomega}{\intttau\!\intomega}
\newcommand{\intoTomega}{\intoT\!\intomega}
\newcommand{\intinfomega}{\int_0^\infty\!\!\intomega}

\newcommand{\intromega}{\int_{\romega}\!} 

\newcommand{\Lo}[1][1]{L^{#1}(\Omega)} 
\newcommand{\W}[1][1,2]{W^{#1}(\Omega)}

\newcommand{\LSp}[2]{L^{#1\;\!}\!\left(#2\right)} 
\newcommand{\LSpn}[2]{L^{#1\;\!}\!(#2)}
\newcommand{\LSpb}[2]{L^{#1\;\!}\!\big(#2\big)}
\newcommand{\LSploc}[2]{L_{loc}^{#1}\!\left(#2\right)} 

\newcommand{\LSplocb}[2]{L_{loc}^{#1}\big(#2\big)}

\newcommand{\WSpn}[2]{W^{#1\,}\!(#2)}
\newcommand{\CSp}[2]{C^{#1}\!\left(#2\right)}

\newcommand{\CSpnl}[2]{C^{#1}\!\,(#2)} 

\newcommand{\R}{\mathbb{R}}
\newcommand{\N}{\mathbb{N}}
\newcommand{\nfrac}[2]{{\nicefrac{#1}{#2}}}
\newcommand{\nepkap}{n_\epsi^{(\kappa)}}
\newcommand{\cepkap}{c_\epsi^{(\kappa)}}
\newcommand{\uepkap}{u_\epsi^{(\kappa)}}
\newcommand{\nkap}{n^{(\kappa)}}
\newcommand{\ckap}{c^{(\kappa)}}
\newcommand{\ukap}{u^{(\kappa)}}

\newcommand{\Td}{{T_{\diamond}}}

\author{Tobias Black\thanks{Institut f\"ur Mathematik, Universit\"at Paderborn, Warburger Str. 100, 33098 Paderborn, Germany; email: \mbox{tblack@math.upb.de}}}
\title{The Stokes limit in a three-dimensional chemotaxis-Navier-Stokes system}
\date{}

\begin{document}
\maketitle
\begin{abstract}
\noindent
{\textbf{Abstract:} We consider initial-boundary value problems for the $\kappa$-dependent family of chemotaxis-(Navier--)Stokes systems
\begin{align*}
\left\{
\begin{array}{r@{\,}c@{\,}c@{\ }l@{\quad}l@{\quad}l@{\,}c}
n_{t}&+&u\cdot\!\nabla n&=\Delta n-\nabla\!\cdot(n\nabla c),\ &x\in\Omega,& t>0,\\
c_{t}&+&u\cdot\!\nabla c&=\Delta c-cn,\ &x\in\Omega,& t>0,\\
u_{t}&+&\kappa(u\cdot\nabla)u&=\Delta u+\nabla P+n\nabla\phi,\ &x\in\Omega,& t>0,\\
&&\nabla\cdot u&=0,\ &x\in\Omega,& t>0,
\end{array}\right.
\end{align*}
in a bounded domain $\Omega\subset\mathbb{R}^3$ with smooth boundary and given potential function $\phi\in\CSp{1+\beta}{\bomega}$ for some $\beta>0$. It is known that for fixed $\kappa\in\R$ an associated initial-boundary value problem possesses at least one global weak solution $(n^{(\kappa)},c^{(\kappa)},u^{(\kappa)})$, which after some waiting time becomes a classical solution of the system. In this work we will show that upon letting $\kappa\to0$ the solutions $(n^{(\kappa)},c^{(\kappa)},u^{(\kappa)})$ converge towards a weak solution of the Stokes variant $(\kappa=0)$ of the systems above with respect to the strong topology in certain Lebesgue and Sobolev spaces. 

\noindent We thereby extend the recently obtained result on the Stokes limit process for classical solutions in the two-dimensional setting to the more intricate three-dimensional case.
\noindent
}\\[0.1cm]


\end{abstract}


\newpage
\section{Introduction}\label{sec1:intro}
The migration towards nutrients is a driving force of nature and even some of the smallest of organisms try to move to better environmental conditions indicated by an increase in concentration of an attracting chemical substance. This phenomenon of biased movement  along a chemical signal gradient is known as chemotaxis and can be observed for a wide array of aerobic bacteria such as \emph{Bacillus subtilis}. Experiments on colonies of \emph{Bacillus subtilis} suspended in a sessile drop of water undertaken in (\cite{tuval2005bacterial}) showed the emergence of plume-like structures and large-scale convection patterns. For theoretical descriptions of the processes involved the authors of said study proposed an extension for the classical Keller--Segel chemotaxis model (\cite{KS70}) capturing the feedback between liquid environment and bacteria, which, upon prototypical choices for the system parameters, can be expressed as
\begin{align*}
\left\{
\begin{array}{r@{\,}c@{\,}c@{\ }l@{\quad}l@{\quad}l@{\,}c}
n_{t}&+&u\cdot\!\nabla n&=\nabla\!\cdot\big(D(n)\nabla n-n\nabla c\big),\ &x\in\Omega,& t>0,\\
c_{t}&+&u\cdot\!\nabla c&=\Delta c-cn,\ &x\in\Omega,& t>0,\\
u_{t}&+&\kappa(u\cdot\nabla)u&=\Delta u+\nabla P+n\nabla\phi,\ &x\in\Omega,& t>0,\\
&&\nabla\cdot u&=0,\ &x\in\Omega,& t>0.
\end{array}\right.
\end{align*}
Herein, the unknown functions $n$ and $c$ represent the density of bacteria and the concentration of the attracting chemical, respectively, $u$ denotes the fluid velocity field and $P$ symbolizes the pressure of the fluid. An archetypical choice for the diffusion rate $D(n)$ is $D\equiv const.$ and the function $\phi$ describes a given gravitational potential, capturing the effect that spots with a high density of bacteria in the fluid are heavier than ones with a low density and tend to sink down.

\noindent{\textbf{Neglecting the fluid convection term.} 
The interplay of the chemotaxis- and Navier--Stokes-equations present in the model poses a very challenging mathematical problem. In particular for $\Omega\subset\R^3$ neither of them is understood completely. For instance, working in the fluid-free three-dimensional setting, obtained upon letting $u\equiv0$ in the system above, global bounded classical solutions were only obtained under the assumption that the initial chemical concentration $\|c(\cdot,0)\|_{\Lo[\infty]}$ is small (\cite{Tao-oxygencons-JMAA11}). In contrast, for arbitrary initial data global weak solutions have been shown to exist, which become smooth and classical after some waiting time (\cite{TaoWin12_evsmooth}).
On the other hand, existence theory for the Navier--Stokes equations, which has been garnering lots of interest for the better part of a century, beyond mere global weak solutions also remains dependent on various assumptions in the three-dimensional setting (\cite{sohr}). Correspondingly, the known results for the chemotaxis-Navier--Stokes systems for arbitrary initial data also only cover global existence of weak solutions (\cite{win_globweak3d-AHPN16}) and eventual smoothing properties (\cite{win_chemonavstokesfinal_TransAm17}). Even in more favorable scenarios, where the diffusion process is enhanced at large cell densities as e.g. incorporated by the choice $D(s)=s^{m-1}$, $s>0$, with $m>1$, only weak solutions could be established, as indicated by the global existence results of \cite{Mizukami-CTNSnonlindiff-2018,ZhangLi-GlobWeakSol-JDE15} covering $m>\frac{2}{3}$.

Accordingly, a wide array of studies dedicated to the mathematical analysis of chemotaxis-fluid interaction mainly concentrates on systems where the fluid evolution is described by the Stokes equation obtained by letting $\kappa=0$, i.e. 
\begin{align}\label{CTStokes}\tag{$\Lambda_0$}
\left\{
\begin{array}{r@{\,}c@{\,}r@{\ }l@{\quad}l@{\quad}l@{\,}c}
n_{t}&+&u\cdot\!\nabla n&=\Delta n-\nabla\!\cdot(n\nabla c),\ &x\in\Omega,& t>0,\\
c_{t}&+&u\cdot\!\nabla c&=\Delta c-cn,\ &x\in\Omega,& t>0,\\
&&u_{t}&=\Delta u+\nabla P+n\nabla\phi,\ &x\in\Omega,& t>0,\\
&&\nabla\cdot u&=0,\ &x\in\Omega,& t>0.
\end{array}\right.
\end{align}
In this setting substantially stronger results besides mere global existence (\cite{DuanXiang-NoteOn-IMRN14,Win-glob-large-data_CPDE12}) have be shown (see e.g \cite{Chae2014,Duan2010,kozono15} and \cite[Section 4.1]{BBWT15} for an additional non-exhaustive overview). The reasoning behind the neglection of the convection term, however, mostly originates from experimental observations indicating Reynolds numbers of order $\mathcal{R}\approx 10^{-4}$ (\cite{Mendelson1999}) for the bacteria in question. Rigorous mathematical results appear to be mostly lacking. In fact, only recently it was shown in the two-dimensional setting that upon taking $\kappa\to0$ the global classical solution $\big(\nkap,\ckap,\ukap\big)$ of the chemotaxis-Navier--Stokes system convergences uniformly in time towards the global classical solution $(n^{(0)},c^{(0)},u^{(0)}\big)$ of \eqref{CTStokes} in the sense that there exist $C>0$ and $\mu>0$ such that whenever $\kappa\in(-1,1)$,
\begin{align*}
\big\|\nkap(\cdot,t)-n^{(0)}(\cdot,t)\big\|_{\Lo[\infty]}+\big\|\ckap(\cdot,t)-c^{(0)}(\cdot,t)\big\|_{\Lo[\infty]}+\big\|\ukap(\cdot,t)-u^{(0)}(\cdot,t)\big\|_{\Lo[\infty]}\leq C|\kappa|e^{-\mu t}
\end{align*}
holds for all $t>0$ (\cite{WangWinklerXiang-smallconvection-MathZ18}).

\noindent{\textbf{Main results.} Motivated by the temporally uniform convergence result for the limit $\kappa\to0$ from \cite{WangWinklerXiang-smallconvection-MathZ18} we aspire to quantify the effect of the Stokes approximation in the more intricate three dimensional setting beyond the expected mere time-local convergence. Before we take a brief look at the major challenges entailed by the increased space dimension, let us specify the framework and the main result obtained in this work. Under the assumptions that $\Omega\subset\R^3$ is a bounded domain with smooth boundary and that $\kappa\in[-1,1]$ we will consider 
\begin{align}\label{CTNScons}\tag{$\Lambda_\kappa$}
\left\{
\begin{array}{r@{\,}c@{\,}c@{\ }l@{\quad}l@{\quad}l@{\,}c}
n_{t}&+&u\cdot\!\nabla n&=\Delta n-\nabla\!\cdot(n\nabla c),\ &x\in\Omega,& t>0,\\
c_{t}&+&u\cdot\!\nabla c&=\Delta c-cn,\ &x\in\Omega,& t>0,\\
u_{t}&+&\kappa(u\cdot\nabla)u&=\Delta u+\nabla P+n\nabla\phi,\ &x\in\Omega,& t>0,\\
&&\nabla\cdot u&=0,\ &x\in\Omega,& t>0,
\end{array}\right.
\end{align}
with boundary conditions
\begin{align}\label{BC}
\nabla n(x,t)\cdot\nu=0,\quad\nabla c(x,t)\cdot\nu=0\quad\text{and}\quad u(x,t)=0\qquad\text{for }x\in\romega\text{ and }t>0,
\end{align}
and initial conditions
\begin{align}\label{IC}
n(x,0)=n_0(x),\quad c(x,0)=c_0(x),\quad u(x,0)=u_0(x),\quad x\in\Omega,
\end{align}
where
\begin{align}\label{phidef}
\phi\in C^{1+\beta}(\bomega)\quad\text{for some }\beta>0.
\end{align}
Moreover, we assume the initial data to satisfy
\begin{align}\label{IR}
\left\{\begin{array}{r@{\,}l}
n_0&\in\CSp{0}{\bomega}\quad\text{is nonnegative with } n_0\not\equiv0,\\
c_0&\in\W[1,\infty]\quad\text{ with }c_0>0\text{ in }\bomega,\\
u_0&\in D(A^\alpha)\quad\text{for some } \alpha\in(\tfrac{3}{4},1),
\end{array}\right.
\end{align}
where $A:=-\mathcal{P}\Delta$ denotes the realization of the Stokes operator in $\LSpn{2}{\Omega;\R^3}$ under homogeneous Dirichlet boundary conditions with its domain $D(A):=\WSpn{2,2}{\Omega;\R^3}\cap W_0^{1,2}(\Omega;\R^3)\cap L_\sigma^{2}\!\left(\Omega\right)$. Herein, $L_\sigma^{2}\!\left(\Omega\right):=\left\{\varphi\in\LSpn{2}{\Omega;\R^3}\,\vert\,\nabla\cdot\varphi=0\right\}$ represents the Hilbert space of solenoidal vector fields in $\Lo[2]$ and $\mathcal{P}$ stands for the Helmholtz projection of $\LSpn{2}{\Omega;\R^3}$ onto $L_\sigma^{2}\left(\Omega\right)$. Accordingly, we also abbreviate $W_{0,\sigma}^{1,p}\left(\Omega\right):=W_0^{1,p}\!(\Omega;\R^3)\cap L_\sigma^{2}\!\left(\Omega\right)$ and $C_{0,\sigma}^\infty(\Omega):=C_0^{\infty}(\Omega;\R^3)\cap L_\sigma^{2}\left(\Omega\right)$.

With the framework and notations clarified, we can now precisely state our main result.

\begin{theorem}\label{theo:1}
Let $\Omega\subset\R^3$ be a bounded and smooth domain and suppose that $\phi$ and $n_0$, $c_0$, $u_0$ comply with \eqref{phidef} and \eqref{IR}, respectively.
Let \begin{align*}
X:=&\ \LSpb{\infty}{(0,\infty);\Lo[1]}\cap\LSplocb{\frac53}{\bomega\times[0,\infty)}\cap\LSplocb{\frac54}{[0,\infty);\W[1,\frac54]}\\&\ \qquad\times\LSpb{\infty}{\Omega\times(0,\infty)}\cap\LSploc{4}{[0,\infty);\W[1,4]}\\&\ \qquad\qquad\times\LSploc{\infty}{[0,\infty);L_\sigma^2(\Omega)}\cap\LSplocb{2}{[0,\infty);W_{0,\sigma}^{1,2}(\Omega)}.
\end{align*}
Then there exist a family 
$\big\{\big(\nkap,\ckap,\ukap\big)\big\}_{\kappa\in[-1,1]}\subset X$ of global weak solutions, in the sense of Definition \ref{def:sol} below, to the corresponding family of chemotaxis-Navier--Stokes systems \eqref{CTNScons},\eqref{BC},\eqref{IC} and $\Td>0$ such that $(\nkap,\ckap,\ukap)$ together with some $P^{(\kappa)}\in\CSp{1,0}{\bomega\times(\Td,\infty)}$ solve \eqref{CTNScons},\eqref{BC},\eqref{IC} classically in $\Omega\times(\Td,\infty)$. Moreover, for any null sequence $(\kappa_j)_{j\in\N}\subset[-1,1]$ there exist a subsequence $(\kappa_{j_k})_{k\in\N}$ and a global weak solution $(n,c,u)\in X$ of the chemotaxis-Stokes system \eqref{CTStokes},\eqref{BC},\eqref{IC}, such that
\begin{align}\label{eq:thm1}
(n^{(\kappa_{j_k})}-n) &\to 0\quad\text{in }\LSpb{p_1}{\Omega\times(0,\infty)}\quad\text{for any }p_1\in[1,\tfrac53),\nonumber\\
(\nabla n^{(\kappa_{j_k})}-\nabla n) &\to 0\quad\text{in }\LSpb{p_2}{\Omega\times(0,\infty);\R^3}\quad\text{for any }p_2\in[1,\tfrac54),\nonumber\\
(c^{(\kappa_{j_k})}-c) &\to 0\quad\text{in }\LSpb{q_1}{\Omega\times(0,\infty)}\quad\text{for any }q_1\in[1,\infty),\\
(\nabla c^{(\kappa_{j_k})}-\nabla c) &\to 0\quad\text{in }\LSpb{q_2}{\Omega\times(0,\infty);\R^3}\quad\text{for any }q_2\in[1,4),\nonumber\\ 
(u^{(\kappa_{j_k})}-u) &\to 0\quad\text{in }\LSpb{r_1}{\Omega\times(0,\infty);\R^3}\quad\text{for any }r_1\in[1,\tfrac{10}{3}),\nonumber\\
(\nabla u^{(\kappa_{j_k})}-\nabla u) &\to 0\quad\text{in }\LSpb{r_2}{\Omega\times(0,\infty);\R^{3\times 3}}\quad\text{for any }r_2\in[1,2)\nonumber
\end{align}
as $\kappa_{j_k}\to0$, and such that $(n,c,u)$ together with some $P\in\CSp{1,0}{\bomega\times(\Td,\infty)}$ solve \eqref{CTStokes},\eqref{BC},\eqref{IC} classically in $\Omega\times(\Td,\infty)$.
\end{theorem}
\noindent{\textbf{Mathematical challenges and the approach.}
In the two-dimensional setting investigated in \cite{WangWinklerXiang-smallconvection-MathZ18}, it is known that \eqref{CTNScons} already admits a classical solution on $\Omega\times(0,\infty)$, which in turn allows for testing procedures immediately targeting the quasi-energy functional
\begin{align*}
\intomega \nkap\ln \nkap+\frac{1}{2}\intomega\frac{\big|\nabla \ckap\big|^2}{\ckap}+\eta\intomega\big|\ukap\big|^2
\end{align*}
for large $\eta>0$ independent of $\kappa\in[-1,1]$ to derive, after a some bootstrapping, $\kappa$-independent bounds in $\CSp{1}{\bomega}\times\CSp{2}{\bomega}\times D(A^\alpha)$ uniform in time. These bounds, when combined with decay properties of \eqref{CTNScons}, then become the driving force of the exponential stabilization featured in \cite{WangWinklerXiang-smallconvection-MathZ18}. In stark contrast, in the current three-dimensional framework we cannot utilize a corresponding quasi-energy functional immediately, as for \eqref{CTNScons} only the global existence of a weak solution obtained by a limiting procedure from approximating systems is known (\cite{win_globweak3d-AHPN16}). To transfer any reasonable information to this weak solution, however, we have to ensure that the precompactness properties used in the limit procedure are independent of $\kappa$. Even though the methods behind the derivation of the corresponding bounds are known (the same quasi-energy as above is exploited for the approximate system), their possible dependence on $\kappa$ has not yet been ruled out and will be inspected in Sections \ref{sec2:locex} and \ref{sec3:epsilimit}. While the strong convergence properties entailed by these bounds (due to the independence of $\kappa$) would also entail a time-local convergence in certain $L^p$ spaces in the limit $\kappa\to0$, we strive for a stronger convergence result global in time. To expand the knowledge, however, we will need to meticulously adjust the analytic machinery behind the eventual smoothness results of \cite{win_chemonavstokesfinal_TransAm17,Lan-Longterm_M3AS16} in order to be able to carefully track the possible $\kappa$-dependence in the eventual smallness of $\cepkap$, the eventual regularity estimates for $\nepkap$ and $\uepkap$ and their eventual stabilization properties presented in Sections \ref{sec4:decay} -- \ref{sec6:ev-dec-n-u}. We can then utilize maximal Sobolev regularity estimates for the Stokes and Neumann heat-semigroups to obtain an eventual smoothing time $\Td>0$, which does not depend on $\kappa$, ensuring that the triple $\big(\nkap,\ckap,\ukap\big)$, obtained in the limit $\epsi\to0$, solves \eqref{CTNScons} classically on $\Omega\times(\Td,\infty)$ (Section \ref{sec7:evsmooth}). Section \ref{sec8:expdec} will then be devoted to gain insight in exponential decay estimates valid starting from the smoothing time $\Td>0$ and finally in Section 9, we will take $\kappa\to0$ to obtain Theorem \ref{theo:1}.
\setcounter{equation}{0} 
\section{Preliminaries. Weak solutions and a priori information for a family of approximating systems}\label{sec2:locex}

Before we start with our detailed analysis let us also briefly specify what constitutes a weak solution as mentioned in Theorem \ref{theo:1}. In the following definition, adapted from \cite{win_globweak3d-AHPN16}, we merely prescribe the weakest regularity necessary to ensure that all integrals in the equalities below are well defined. The solutions constructed later, however, will satisfy considerably stronger regularity assumptions. 
\begin{definition}\label{def:sol}
For $\kappa\in[-1,1]$ a triple $(\nkap,\ckap,\ukap)$ of functions
\begin{align*}
\nkap\in\LSploc{1}{[0,\infty);\W[1,1]},\quad
\ckap\in\LSploc{1}{[0,\infty);\W[1,1]},\quad
\ukap\in\LSplocb{1}{[0,\infty);W_0^{1,1}(\Omega;\R^3)},
\end{align*} 
satisfying $\nkap\geq0$, $\ckap\geq0$ and $\nabla\cdot\ukap=0$ a.e. in $\bomega\times[0,\infty)$, $\nkap\ckap\in\LSploc{1}{\bomega\times[0,\infty)}$, and $\kappa\ukap\otimes\ukap\in\LSplocb{1}{\bomega\times[0,\infty);\R^{3\times3}}$ with
\begin{align*}
\nkap\nabla\ckap,\quad \text{and}\quad \nkap\ukap\quad \text{as well as}\quad \ckap\ukap\quad \text{belonging to}\quad\LSplocb{1}{\bomega\times[0,\infty);\R^3},
\end{align*}
will be called a weak solution of the system \eqref{CTNScons}, \eqref{BC} and \eqref{IC}, if the equality
\begin{align*}
-\intinfomega\nkap\varphi-\intomega n_0\varphi(\cdot,0)&=\intinfomega\nkap\ukap\cdot\nabla\varphi-\intinfomega\nabla\nkap\cdot\nabla\varphi+\intinfomega\nkap\nabla\ckap\cdot\nabla\varphi
\intertext{holds for every $\varphi\in C_0^\infty\!\left(\bomega\times[0,\infty)\right)$, if moreover}
-\intinfomega\ckap\psi_t-\intomega c_0\psi(\cdot,0)&=\intinfomega\ckap\ukap\cdot\nabla\psi-\intinfomega\nabla\ckap\cdot\nabla\psi-\intinfomega\nkap\ckap\psi
\intertext{is fulfilled for every $\psi\in C_0^\infty\!\left(\bomega\times[0,\infty)\right)$, and if finally}
-\intinfomega\ukap\cdot\Psi_t-\intomega u_0\cdot\Psi(\cdot,0)&=-\intinfomega\nabla\ukap\cdot\nabla\Psi+\kappa\intinfomega\ukap\otimes\ukap\cdot\nabla\Psi+\intinfomega\nkap\Psi\cdot\nabla\phi
\end{align*}
is valid for every $\Psi\in C_{0}^\infty\!\left(\Omega\times[0,\infty);\R^3\right)$ satisfying $\nabla\cdot\Psi\equiv0$.
\end{definition}

Weak solutions to \eqref{CTNScons}, in the sense above, will be constructed as limit objects from a family of appropriately regularized systems. The regularization we incorporate for our problem has previously (and in a more general fashion) been employed in \cite{Lan-Longterm_M3AS16,win_globweak3d-AHPN16,win_chemonavstokesfinal_TransAm17}. To be precise, for $\epsi\in(0,1)$ and $\kappa\in[-1,1]$ we will consider
\begin{align}\label{approxprop}\tag{$\Lambda_{\epsi,\kappa}$}
\left\{
\begin{array}{r@{\,}l@{\quad}l@{\quad}l@{\,}c}
n_{\epsi t}^{(\kappa)}+\uepkap\cdot\!\nabla \nepkap&=\Delta \nepkap-\nabla\!\cdot\Big(\frac{\nepkap}{1+\epsi\nepkap}\nabla \cepkap\Big),\ &x\in\Omega,& t>0,\\
c_{\epsi t}^{(\kappa)}+\uepkap\cdot\!\nabla \cepkap&=\Delta \cepkap-\frac{1}{\epsi}\ln\big(1+\epsi\nepkap\big)\cepkap,\ &x\in\Omega,& t>0,\\
u_{\epsi t}^{(\kappa)}+\kappa(Y_\epsi\uepkap\cdot\nabla)\uepkap&=\Delta \uepkap+\nabla P_\epsi^{(\kappa)}+\nepkap\nabla\phi,\ &x\in\Omega,& t>0,\\
\nabla\cdot \uepkap&=0,\ &x\in\Omega,& t>0,\\
\partial_\nu \nepkap=0,&\ \!\!\;\qquad \partial_\nu \cepkap=0,\qquad\qquad\quad \uepkap=0,\ &x\in\romega,& t>0,\\
\nepkap(x,0)=n_0(x),&\quad \cepkap(x,0)=c_0(x),\quad \uepkap(x,0)=u_0(x),\ &x\in\Omega,&
\end{array}\right.
\end{align}
where for $\epsi\in(0,1)$ $Y_\epsi$ denotes the Yosida approximation (\cite{MiyakawaSohr-MathZ88,sohr}) given by
\begin{align*}
Y_\epsi\varphi:=(1+\epsi A)^{-1}\varphi\quad\text{for }\varphi\in L_\sigma^2(\Omega).
\end{align*}
Let us also note that
\begin{align}\label{eq:approxest}
\frac{1}{2}\min\{s,1\}\leq \frac{1}{\epsi}\ln(1+\epsi s)\leq s\quad\text{for all }s\geq0\text{ and all }\epsi\in(0,1),
\end{align}
which, due to nonnegativity of $\nepkap$ we will establish later, are two useful estimates for one of the terms appearing in the second equation of \eqref{approxprop}, which we will use on multiple occasions throughout the paper.

Now, let us start our analysis by gathering basic results for the family of approximating systems, most of which has already been discussed in works with fixed $\kappa=1$ and can be obtained in well-known manner. Nevertheless, we have to ascertain that all of these familiar properties are $\kappa$-independent and therefore will take a closer look at some (parts) of the proofs involved.

\begin{lemma}
\label{lem:locex}
Let $q>3$. For any $\epsi\in(0,1)$ and $\kappa\in[-1,1]$ there exist $\Tmek\in(0,\infty]$ and a unique triplet $(\nepkap,\cepkap,\uepkap)$ of functions satisfying
\begin{align*}
\nepkap&\in\CSpnl{0}{\bomega\times[0,\Tmek)}\cap \CSpnl{2,1}{\bomega\times(0,\Tmek)},\\\cepkap&\in\CSpnl{0}{\bomega\times[0,\Tmek)}\cap\CSpnl{2,1}{\bomega\times(0,\Tmek) }\cap\LSpb{\infty}{(0,\Tmek);\W[1,q]},\\\uepkap&\in\CSpnl{0}{\bomega\times[0,\Tmek);\R^3}\cap \CSpnl{2,1}{\bomega\times(0,\Tmek);\R^3},
\end{align*}
which together with some $P_\epsi^{(\kappa)}\in\CSpnl{1,0}{\Omega\times(0,\Tmek)}$ solves \eqref{approxprop} classically in $\Omega\times(0,\Tmek)$. In addition,
\begin{align*}
\mbox{if }\Tmek<\infty,\mbox{ then }\|\nepkap(\cdot,t)\|_{\Lo[\infty]}+\|\cepkap(\cdot,t)\|_{\W[1,q]}+\|&A^\alpha\uepkap(\cdot,t)\|_{\Lo[2]}\to\infty\mbox{ as }t\nearrow\Tmek\\
&\qquad\mbox{ for all }\alpha\in(\tfrac{3}{4},1).
\end{align*}
 The triplet $\big(\nepkap,\cepkap,\uepkap)$ moreover satisfies $\nepkap\geq 0$ and $\cepkap>0$ in $\bomega\times[0,\Tmek)$, as well as
\begin{align}\label{eq:locex-bounds}
\intomega \nepkap(\cdot,t)=\intomega n_0\quad\text{and}\quad\|\cepkap(\cdot,t)\|_{\Lo[\infty]}\leq\|c_0\|_{\Lo[\infty]}\quad\text{for all }t\in[0,\Tmek),
\end{align}
and the mapping $t\mapsto\big\|\cepkap(\cdot,t)\big\|_{\Lo[\infty]}$ is nonincreasing on $(0,\infty)$.
\end{lemma}

\begin{bew}
The proof draws on a standard reasoning involving semigroup estimates, Banach's fixed point theorem employed to a closed subset of $\LSp{\infty}{(0,T);\CSp{0}{\bomega}\times\W[1,q]\times D(A^\alpha)}$ and parabolic regularity theory. We refer the reader to \cite[Lemma 2.1]{Win-glob-large-data_CPDE12} for a detailed proof of the existence of a unique local solution, the extensibility criterion and the nonnegativity and positivity properties in a closely related setting. The conservation of mass $\intomega \nepkap=\intomega n_0$ on $(0,\Tmek)$ then follows directly from integrating the first equation of \eqref{approxprop}, whereas the nonincreasing property of $t\mapsto\big\|\cepkap(\cdot,t)\big\|_{\Lo[\infty]}$ on $(0,\infty)$ and bound for $\big\|\cepkap\big\|_{\Lo[\infty]}$ are an immediate consequence of the parabolic comparison principle employed to the second equation of \eqref{approxprop}.
\end{bew}

Since $\kappa$ only impacts the third equation of \eqref{approxprop} directly, we can, without any necessary change, adopt the results from \cite[Lemmas 2.6 and 2.8]{Lan-Longterm_M3AS16} and \cite[Lemma 3.4]{win_globweak3d-AHPN16} to obtain the following.

\begin{lemma}
\label{lem:diff-ineq-1}
There exists $K_0>0$ such that for all $\epsi\in(0,1)$ and all $\kappa\in[-1,1]$ the solution $(\nepkap,\cepkap,\uepkap)$ of \eqref{approxprop} satisfies
\begin{align*}
\frac{\intd}{\intd t}&\left(\intomega\nepkap\ln\nepkap+\frac{1}{2}\intomega\frac{\big|\nabla\cepkap\big|^2}{\cepkap}\right)\\&\hspace*{2.38cm}+\frac{1}{K_0}\left(\intomega\frac{\big|\nabla\nepkap\big|^2}{\nepkap}+\intomega\frac{\big|D^2\cepkap\big|^2}{\cepkap}+\intomega\frac{\big|\nabla\cepkap\big|^4}{{\cepkap}^3}\right)\leq K_0\intomega\big|\nabla\uepkap\big|^2+K_0
\end{align*}
on $(0,\Tmek)$.
\end{lemma}

\begin{bew}
Since the well-established testing procedures used to derive this inequality do not depend on $\kappa$ in any way, we refer the reader to the detailed proofs in \cite[Lemmas 2.6 and 2.8]{Lan-Longterm_M3AS16} (with $\kappa=1$) and \cite[Lemma 3.4]{win_globweak3d-AHPN16} (in convex domains with $\kappa=1$).
\end{bew}

Moreover, due to $\uepkap$ being divergence free, testing the third equation against $\uepkap$ itself also removes any dependence on $\kappa$ and hence we readily transfer the result from \cite[Lemma 2.9]{Lan-Longterm_M3AS16} to our setting.

\begin{lemma}
\label{lem:diff-ineq-2}
For any $\epsi\in(0,1)$ and $\kappa\in[-1,1]$ the solution $(\nepkap,\cepkap,\uepkap)$ of \eqref{approxprop} satisfies
\begin{align*}
\frac12\frac{\intd}{\intd t}\intomega\big|\uepkap\big|^2+\intomega\big|\nabla\uepkap\big|^2=
\intomega\nepkap\nabla\phi\cdot\uepkap
\end{align*}
on $(0,\Tmek)$.
\end{lemma}

\begin{bew}
Since $\nabla\cdot\uepkap=0$ on $\Omega\times(0,\Tmek)$ also implies that $\nabla\cdot Y_\epsi\uepkap=0$ on $\Omega\times(0,\Tmek)$, we have
\begin{align*}
\kappa\intomega\big(Y_\epsi\uepkap\cdot\nabla\big)\uepkap\cdot\uepkap=-\kappa\intomega\nabla\cdot\big(Y_\epsi\uepkap\big)\big|\uepkap\big|^2-\frac{\kappa}{2}\intomega Y_\epsi\uepkap\cdot\nabla\big|\uepkap\big|^2
=0
\end{align*}
on $\Omega\times(0,\Tmek)$. Thus, we find that by multiplying the third equation of \eqref{approxprop} by $\uepkap$ and integrating by parts
\begin{align*}
\frac{1}{2}\frac{\intd}{\intd t}\intomega\big|\uepkap\big|^2+\intomega\big|\nabla\uepkap\big|^2=\intomega\nepkap\nabla\phi\cdot\uepkap
\end{align*}
is valid on $(0,\Tmek)$.
\end{bew}

Combination of the previous two lemmas now yields uniform a priori estimates which will be the basis for the remainder of our regularity analysis.

\begin{lemma}
\label{lem:bounds}
Let $K_0>0$ be provided by Lemma \ref{lem:diff-ineq-2}. There exists $K_1>0$ such that for all $\epsi\in(0,1)$ and each $\kappa\in[-1,1]$ the solution $(\nepkap,\cepkap,\uepkap)$ of \eqref{approxprop} satisfies
\begin{align*}
\intomega\nepkap\ln\nepkap+\frac{1}{2}\intomega\frac{\big|\nabla \cepkap\big|^2}{\cepkap}+K_0\intomega\big|\uepkap\big|^2\leq K_1
\end{align*}
on $(0,\Tmek)$ and
\begin{align*}
\inttauomega\frac{\big|\nabla\nepkap\big|^2}{\nepkap}+\inttauomega\frac{\big|D^2\cepkap\big|^2}{\cepkap}+\inttauomega\frac{\big|\nabla\cepkap\big|^4}{{\cepkap}^3}+\inttauomega\big|\nabla\uepkap\big|^2+\inttauomega\big|\nabla\cepkap\big|^4\leq K_1
\end{align*}
for all $t\in(0,\Tmek-\tau)$, where $\tau:=\min\big\{1,\tfrac12\Tmek\big\}$.
\end{lemma}

\begin{bew}
(Compare \cite[Lemmas 2.10 and 2.11]{Lan-Longterm_M3AS16} and \cite[Lemmas 3.6 and 3.8]{win_globweak3d-AHPN16}.)
Adding up suitable multiples of the differential inequalities from Lemma \ref{lem:diff-ineq-1} and Lemma \ref{lem:diff-ineq-2} we find that for any $\epsi\in(0,1)$, $\kappa\in[-1,1]$
\begin{align}\label{eq:bounds-eq1}
\frac{\intd}{\intd t}\bigg(\intomega\nepkap&\ln\nepkap+\frac{1}{2}\intomega\frac{\big|\nabla\cepkap\big|^2}{\cepkap}+K_0\intomega\big|\uepkap\big|^2\bigg)+K_0\intomega\big|\nabla\uepkap\big|^2\nonumber\\
&+\frac{1}{K_0}\bigg(\intomega\frac{\big|\nabla\nepkap\big|^2}{\nepkap}+\intomega\frac{\big|D^2\cepkap\big|^2}{\cepkap}+\intomega\frac{\big|\nabla\cepkap\big|^4}{{\cepkap}^3}\bigg)\leq \intomega\nepkap\nabla\phi\cdot\uepkap+K_0
\end{align}
holds on $(0,\Tmek)$. To estimate the right-hand side further, we make use of the boundedness of $\nabla\phi$ and Hölder's inequality, the embedding $W^{1,2}_0(\Omega)\hookrightarrow\Lo[6]$ and the Poincaré inequality to obtain $C_1>0$ such that for each $\epsi\in(0,1)$, $\kappa\in[-1,1]$ and all $t\in(0,\Tmek)$ we have
\begin{align*}
\intomega\nepkap\nabla\phi\cdot\uepkap\leq\big\|\nabla\phi\big\|_{\Lo[\infty]}\big\|\nepkap\big\|_{\Lo[\frac65]}\big\|\uepkap\big\|_{\Lo[6]}\leq C_1\big\|\nabla\phi\big\|_{\Lo[\infty]}\big\|\nepkap\big\|_{\Lo[\frac65]}\big\|\nabla\uepkap\big\|_{\Lo[2]}.
\end{align*}
Here, we employ Young's inequality to find that
\begin{align}\label{eq:bounds-eq2}
\intomega\nepkap\nabla\phi\cdot\uepkap\leq \frac{K_0}{2}\big\|\nabla\uepkap\big\|_{\Lo[2]}^2+\frac{C_1^2\|\nabla\phi\|_{\Lo[\infty]}^2}{2 K_0}\big\|\nepkap\big\|_{\Lo[\frac{6}{5}]}^2
\end{align}
on $(0,\Tmek)$. According to the \GNI\ there is some $C_2>0$ such that
\begin{align*}
\|\varphi\|_{\Lo[\frac{12}{5}]}^4\leq C_2\|\nabla\varphi\|_{\Lo[2]}\|\varphi\|_{\Lo[2]}^3+C_2\|\varphi\|_{\Lo[2]}^4
\end{align*}
holds for all $\varphi\in\W[1,2]$ and hence, in light of the mass conservation $\intomega\nepkap=\intomega n_0$ for any $\epsi\in(0,1)$, $\kappa\in[-1,1]$ and all $t\in(0,\Tmek)$ from Lemma \ref{lem:locex}, there exists some $C_3>0$ such that for each $\epsi\in(0,1)$ and $\kappa\in[-1,1]$
\begin{align*}
\big\|\nepkap\big\|_{\Lo[\frac{6}{5}]}^2&=\big\|{\nepkap}^\frac{1}{2}\big\|_{\Lo[\frac{12}{5}]}^4\\
&\leq C_2\big\|\nabla{\nepkap}^\frac{1}{2}\big\|_{\Lo[2]}\big\|{\nepkap}^\frac12\big\|_{\Lo[2]}^3+C_2\big\|{\nepkap}^\frac{1}{2}\big\|_{\Lo[2]}^4\leq C_3\big\|\nabla{\nepkap}^{\frac{1}{2}}\big\|_{\Lo[2]}+C_3
\end{align*}
is valid on $(0,\Tmek)$. Employing Young's inequality once more in \eqref{eq:bounds-eq2} we thereby obtain $C_4>0$ such that for any $\epsi\in(0,1)$, $\kappa\in[-1,1]$ and all $t\in(0,\Tmek)$ the inequality
\begin{align*}
\intomega\nepkap\nabla\phi\cdot\uepkap\leq \frac{K_0}{2}\intomega\big|\nabla\uepkap\big|^2+\frac{1}{2K_0}\intomega\frac{\big|\nabla \nepkap\big|^2}{\nepkap}+C_4
\end{align*}
holds. Plugging this into \eqref{eq:bounds-eq1} we find $C_5:=\max\{C_4+K_0,\frac{2}{K_0},2K_0\}>0$ such that for each $\epsi\in(0,1)$, $\kappa\in[-1,1]$ and all $t\in(0,\Tmek)$ the functions
\begin{align*}
y_\epsi^{(\kappa)}(t)&:=\intomega\big(\nepkap\ln\nepkap\big)(\cdot,t)+\frac{1}{2}\intomega\frac{\big|\nabla\cepkap(\cdot,t)\big|^2}{\cepkap(\cdot,t)}+K_0\intomega\big|\uepkap(\cdot,t)\big|^2\\
\mbox{and }h_\epsi^{(\kappa)}(t)&:=\intomega\frac{\big|\nabla\nepkap(\cdot,t)\big|^2}{\nepkap(\cdot,t)}+\intomega\frac{\big|D^2\cepkap(\cdot,t)\big|^2}{\cepkap(\cdot,t)}+\intomega\frac{\big|\nabla\cepkap(\cdot,t)\big|^4}{\big(\cepkap(\cdot,t)\big)^3}+\intomega\big|\nabla\uepkap(\cdot,t)\big|^2
\end{align*} 
satisfy the differential inequality
\begin{align}\label{eq:bounds-eq3}
\frac{\intd}{\intd t}y_\epsi^{(\kappa)}(t)+\frac{1}{C_5}h_\epsi^{(\kappa)}(t)\leq C_5.
\end{align}
Invoking the Poincaré inequality, Young's inequality, the boundedness of $\cepkap$, the inequality $z\ln z\leq\frac{3}{2}z^\frac{5}{4}$ for $z\geq0$, the \GNI\ and the mass conservation of $\nepkap$ from \eqref{eq:locex-bounds} it can be easily checked that there is some $C_6>0$ (independent of $\epsi$ and $\kappa$) such that
\begin{align*}
y_\epsi^{(\kappa)}(t)\leq C_6 h_\epsi^{(\kappa)}(t)+C_6\quad\text{for all }t\in(0,\Tmek).
\end{align*} 
And hence \eqref{eq:bounds-eq3} takes the form
\begin{align*}
\frac{\intd}{\intd t}y_\epsi^{(\kappa)}(t)+\frac{1}{2C_5}h_\epsi^{(\kappa)}(t)+\frac{1}{2C_5C_6}y_\epsi^{(\kappa)}(t)\leq C_5+\frac{1}{2C_5}\quad\text{for all }t\in(0,\Tmek),
\end{align*}
which on the one hand implies for any $\epsi\in(0,1)$, $\kappa\in[-1,1]$ and all $t\in(0,\Tmek)$ that
\begin{align*}
y_\epsi^{(\kappa)}(t)\leq C_7:=\max\Big\{\intomega n_0\ln n_0+\frac12\intomega\frac{|\nabla c_0|^2}{c_0}+K_0\intomega|u_0|^2,2C_5^2C_6+C_6\Big\},
\end{align*}
and, on the other hand, shows upon integration that for each $\epsi\in(0,1)$, $\kappa\in[-1,1]$
\begin{align*}
\frac{1}{2C_5}\int_t^{t+\tau}\!\! h_\epsi^{(\kappa)}(t)\intd t\leq y_\epsi^{(\kappa)}(0)+\big(C_5+\frac{1}{2C_5}\big)\tau\leq C_7+C_5+\frac{1}{2C_5}=:C_8
\end{align*}
is valid for all $t\in(0,\Tmek-\tau)$ with $\tau:=\min\{1,\tfrac12 \Tmek\}$. 
Moreover, drawing on the boundedness of $\cepkap$ obtained in Lemma \ref{lem:locex} we find that
\begin{align*}
\int_t^{t+\tau}\!\!\intomega\big|\nabla\cepkap\big|^4&\leq\sup_{s\in[t,t+\tau]}\big\|\cepkap(\cdot,s)\big\|_{\Lo[\infty]}^3\int_t^{t+\tau}\!\!\intomega\frac{\big|\nabla \cepkap|^4}{{\cepkap}^3}\\&\hspace*{4cm}\leq \|c_0\|^3_{\Lo[\infty]}\int_t^{t+\tau}\!\! h_\epsi^{(\kappa)}(t)\intd t\leq 2C_5C_8\|c_0\|^3_{\Lo[\infty]}
\end{align*}
is valid for each $\epsi\in(0,1)$, $\kappa\in[-1,1]$ and all $t\in(0,\Tmek-\tau)$, completing the proof upon obvious choice of $K_1>0$.
\end{bew}

Assuming a finite maximal existence time, we can now make use of the bounds from the previous lemma to derive a contradiction to the extensibility criterion featured in the local existence result.

\begin{lemma}
\label{lem:globex}
For all $\epsi\in(0,1)$ and $\kappa\in[-1,1]$ the solution to \eqref{approxprop} is global in time, i.e. $\Tmek=\infty$.
\end{lemma}

\begin{bew}
Assuming $\Tmek$ to be finite we will derive a contradiction to the extensibility criterion presented in Lemma \ref{lem:locex}. Reasoning along these lines is common in many related works and can e.g. be found in \cite{win_globweak3d-AHPN16}. For sake of completeness we sketch the main parts of the proof. We first note that, due to $\Tmek<\infty$, Lemma \ref{lem:bounds} provides the existence of $C_1>0$ satisfying
\begin{align}\label{eq:globex1}
\int_0^{\Tmek}\!\!\intomega\big|\nabla\cepkap\big|^4\leq C_1\quad\text{and}\quad\intomega\big|\uepkap(\cdot,t)\big|^2\leq C_1\quad\text{for all }t\in(0,\Tmek).
\end{align}
Testing the first equation of \eqref{approxprop} against $(\nepkap)^3$, we find upon integrating by parts, utilizing the fact that $\frac{s}{1+\epsi s}\leq\frac{1}{\epsi}$ for all $s\geq0$ and invoking Young's inequality that
\begin{align*}
\frac{1}{4}\frac{\intd}{\intd t}\intomega{\nepkap}^4+3\intomega{\nepkap}^2\big|\nabla\nepkap\big|^2&\leq \intomega{\nepkap}^2\big|\nabla\nepkap\big|^2+\intomega\big|\nabla\cepkap\big|^4+\frac{81}{64\epsi^4}\intomega{\nepkap}^4\quad\text{on }(0,\Tmek),
\end{align*}
implying that there is some $C_2>0$ (possibly depending on $\epsi$) such that $\intomega(\nepkap)^4(\cdot,t)\leq C_2$ holds for all $t\in(0,\Tmek)$, according to \eqref{eq:globex1}. Furthermore, in light of the embedding $D(1+\epsi A)=\W[2,2]\cap W_{0,\sigma}^{1,2}\hookrightarrow\Lo[\infty]$ and \eqref{eq:globex1} we obtain $C_3,C_4>0$ satisfying
\begin{align*}
\big\|Y_\epsi\uepkap(\cdot,t)\big\|_{\Lo[\infty]}=\big\|(1+\epsi A)^{-1}\uepkap(\cdot,t)\big\|_{\Lo[\infty]}\leq C_3\big\|\uepkap(\cdot,t)\big\|_{\Lo[2]}\leq C_4\quad\text{for all }t\in(0,\Tmek).
\end{align*}
Hence, testing $u^{(\kappa)}_{\epsi t}+A\uepkap=\mathcal{P}\Big(-\kappa\big(Y_\epsi\uepkap\cdot\nabla\big)\uepkap+\nepkap\nabla\phi\Big)$ against $A\uepkap$ we obtain some $C_5>0$ such that 
\begin{align*}
\frac{1}{2}\frac{\intd}{\intd t}\intomega\big|A^\frac{1}{2}\uepkap\big|^2+\intomega\big|A\uepkap\big|^2\leq\intomega\big|A\uepkap\big|^2+C_5\left(\intomega\big|\nabla\uepkap\big|^2+\intomega {\nepkap}^2\right)\quad\text{on }(0,\Tmek),
\end{align*}
in light of Young's inequality, \eqref{phidef} and the facts that $|\kappa|\leq1$ and $\|\mathcal{P}\varphi\|_{\Lo[2]}\leq \|\varphi\|_{\Lo[2]}$ for all $\varphi\in\Lo[2]$. Since $\intomega|A^\frac{1}{2}\varphi|=\intomega|\nabla\varphi|^2$ for $\varphi\in D(A)$ we thereby find $C_6>0$ fulfilling
\begin{align*}
\intomega\big|\nabla\uepkap(\cdot,t)\big|^2\leq C_6\quad\text{for all }t\in(0,\Tmek).
\end{align*}
Combining these bounds with well-known properties of the Stokes semigroup (see e.g. \cite[p.201]{gig86}) first provides a bound on $\big\|A^\alpha\uepkap(\cdot,t)\big\|_{\Lo[2]}$ for all $t\in(0,\Tmek)$, where $\alpha$ is as in \eqref{IR}. By our choice of $\alpha$, the embedding $D(A^\alpha)\hookrightarrow\Lo[\infty]$ also readily entails an $L^\infty$ bound on the third component. Secondly, combining these bounds with semigroup estimates for the Neumann heat semigroup (e.g. \cite[Lemma 1.3]{win10jde}), \eqref{eq:approxest}, \eqref{eq:locex-bounds} and \eqref{eq:globex1} implies the boundedness of $\big\|\nabla\cepkap(\cdot,t)\big\|_{\Lo[4]}$ for all $t\in(0,\Tmek)$, which upon final combination with Neumann heat semigroup estimates with previous bounds also yields a bound on $\big\|\nepkap(\cdot,t)\big\|_{\Lo[\infty]}$ for all $t\in(0,\Tmek)$, contradicting the extensibility criterion from Lemma \ref{lem:locex}, and hence we conclude $\Tmek=\infty$.
\end{bew}

In a straightforward manner we can also draw on the Gagliardo--Nirenberg and Hölder inequalities to refine the spatio-temporal bounds on the gradient terms in Lemma \ref{lem:bounds} into slightly improved bounds for $\nepkap$, $\nabla\nepkap$ and $\uepkap$. The following lemma will play an important role in deriving the necessary precompactness properties to verify that the objects obtained from the limiting procedure actually constitute a weak solution of our system.

\begin{lemma}
\label{lem:timespace-bounds}
For every $T>0$ there exists $C(T)>0$ such that for any $\epsi\in(0,1)$ and all $\kappa\in[-1,1]$ the solution $(\nepkap,\cepkap,\uepkap)$ of \eqref{approxprop} satisfies
\begin{align*}
\intoTomega{\nepkap}^\frac{5}{3}+\intoTomega\big|\nabla\nepkap\big|^\frac{5}{4}+
\intoTomega\big|\uepkap\big|^\frac{10}{3}\leq C(T).
\end{align*}
\end{lemma}

\begin{bew}
The spatio-temporal bounds follow from immediate applications of the Gagliardo--Nirenberg and Hölder inequalities along with the bounds prepared in Lemma \ref{lem:bounds}. Details on the steps involved are found in \cite[Lemma 3.10]{win_globweak3d-AHPN16}.
\end{bew}

\setcounter{equation}{0} 
\section{Existence of a limit solution family when \texorpdfstring{$\epsi\searrow0$}{epsilon to zero}}\label{sec3:epsilimit}

In preparation of an Aubin--Lions type argument, which is the starting point for our convergence result,
we will require information on the regularity of the time derivatives of our solution components. Again taking care that our estimates do neither depend on $\epsi$ nor on $\kappa$ these bounds on the time derivative will not only be useful for the $\epsi$--limit, but also for the $\kappa$--limit discussed in Section \ref{sec9:kappalimit}.

\begin{lemma}
\label{lem:timereg}
For any $T>0$ there exists $C>0$ such that for each $\epsi\in(0,1)$ and $\kappa\in[-1,1]$ the solution $(\nepkap,\cepkap,\uepkap)$ of \eqref{approxprop} satisfies
\begin{align*}
\intoT\big\|n_{\epsi t}^{(\kappa)}\big\|^{\frac{10}{9}}_{\left(W^{1,10}(\Omega)\right)^*}+
\intoT\Big\|\partial_t\sqrt{\cepkap}\Big\|^{\frac53}_{\left(W^{1,\nfrac52}(\Omega)\right)^*}+
\intoT\big\|u_{\epsi t}^{(\kappa)}\big\|^{\frac54}_{\left(W_{0,\sigma}^{1,5}(\Omega)\right)^*}\leq C.
\end{align*}
\end{lemma}

\begin{bew}
The proof is basically contained in \cite[Lemma 3.11]{win_globweak3d-AHPN16} (where $\kappa=1$ was treated). To ensure that the constant does not depend on $\kappa$, we will illustrate the steps involved for the fluid component. For details regarding the other two estimation procedures (which work along similar lines), we refer the reader to the work mentioned above. Given any fixed $\varphi\in C_{0,\sigma}^\infty(\Omega)$ we test the third equation of \eqref{approxprop} against $\varphi$ and employ Hölder's inequality to obtain that, due to $|\kappa|\leq 1$,
\begin{align*}
&\ \left|\intomega u_{\epsi t}^{(\kappa)}(\cdot,t)\cdot\varphi\right|\\
=\ &\left|-\intomega\nabla\uepkap(\cdot,t)\cdot\nabla\varphi-\kappa\intomega(Y_\epsi\uepkap\otimes\uepkap)(\cdot,t)\cdot\nabla\varphi+\intomega\nepkap(\cdot,t)\nabla\phi\cdot\varphi\right|\\
\leq\ &\left(\big\|\nabla\uepkap(\cdot,t)\big\|_{\Lo[\frac{5}{4}]}+\big\|(Y_\epsi\uepkap\otimes\uepkap)(\cdot,t)\big\|_{\Lo[\frac54]}+\big\|\nepkap(\cdot,t)\nabla\phi\big\|_{\Lo[\frac54]}\right)\|\varphi\|_{\W[1,5]}
\end{align*}
is valid for all $t>0$. In light of \eqref{phidef} we can find $C_1>0$ such that $\|\nabla\phi\|_{\Lo[\infty]}\leq C_1$ and hence Young's inequality entails that, with $C_2:=(1+C_1)>0$, we have
\begin{align}\label{eq:timereg1}
\ &\intoT\big\|u_{\epsi t}^{(\kappa)}(\cdot,t)\big\|^{\frac{5}{4}}_{\left(W_{0,\sigma}^{1,5}(\Omega)\right)^*}\intd t\nonumber\\\leq\ & C_2\intoTomega\big|\nabla\uepkap\big|^\frac{5}{4}+C_2\intoTomega\big|Y_\epsi\uepkap\otimes\uepkap\big|^\frac{5}{4}+C_2\intoTomega{\nepkap}^\frac{5}{4}\nonumber\\
\leq\ & C_2\intoTomega\big|\nabla\uepkap\big|^2+C_2\intoTomega\big|Y_\epsi\uepkap\big|^2+C_2\intoTomega\big|\uepkap\big|^\frac{10}{3}+C_2\intoTomega{\nepkap}^\frac{5}{3}+2C_2|\Omega|T
\end{align}
for all $T>0$. Drawing on the fact that $\|Y_\epsi v\|_{\Lo[2]}\leq\|v\|_{\Lo[2]}$ holds for all $v\in L_\sigma^2(\Omega)$, we may employ Young's inequality once more to estimate $\intoTomega\big|Y_\epsi\uepkap\big|^2\leq\intoTomega\big|\uepkap\big|^\frac{10}{3}+|\Omega|T$ and thus conclude the asserted bound from \eqref{eq:timereg1} in light of Lemmas \ref{lem:bounds} and \ref{lem:timespace-bounds}. 
\end{bew}

With the uniform bounds from Lemma \ref{lem:locex}, Lemma \ref{lem:bounds}, \ref{lem:timespace-bounds} and Lemma \ref{lem:timereg} we are now in the position to obtain limit functions $\nkap$, $\ckap$ and $\ukap$, which fulfill the regularity assumptions and integral equations required to satisfy the weak formulation of \eqref{CTNScons}.

\begin{lemma}
\label{lem:epsi-limit}
There exist a sequence $(\epsi_j)_{j\in\N}\subset(0,1)$ with $\epsi_j\searrow0$ as $j\to\infty$ with the property that for any $\kappa\in[-1,1]$ one can find functions 
\begin{align*}
\nkap&\in\LSploc{\frac53}{\bomega\times[0,\infty)}\quad\text{with}\quad\nabla\nkap\in\LSploc{\frac54}{\bomega\times[0,\infty)},\\
\ckap&\in\LSp{\infty}{\Omega\times(0,\infty)}\quad\text{with}\quad\nabla\ckap\in\LSploc{4}{\bomega\times[0,\infty)},\\
\ukap&\in L^2_{loc}\big([0,\infty);W_{0,\sigma}^{1,2}(\Omega)\big),
\end{align*}
such that the solution $(\nepkap,\cepkap,\uepkap)$ of \eqref{approxprop} satisfies
\begin{alignat}{2}
\nepkap&\to\nkap \qquad&&\text{in }\LSploc{p}{\bomega\times[0,\infty)}\text{ for any }p\in[1,\tfrac53)\text{ and a.e. in }\Omega\times(0,\infty),\label{eq:conv-nepkap}\\
\nabla\nepkap&\wto \nabla\nkap &&\text{in }\LSploc{\frac54}{\bomega\times[0,\infty)},\label{eq:w-conv-nabla-nepkap}\\
\nepkap&\wto\nkap&&\text{in }\LSploc{\frac53}{\bomega\times[0,\infty)},\label{eq:w-conv-nepkap}\\
\cepkap&\to\ckap  &&\text{in }\LSploc{p}{\bomega\times[0,\infty)}\text{ for any }p\in[1,\infty)\text{ and a.e. in }\Omega\times(0,\infty),\label{eq:conv-cepkap}\\
\cepkap&\wsto\ckap&&\text{in }\LSp{\infty}{\Omega\times(0,\infty)},\label{eq:wst-conv-cepkap}\\
\nabla\cepkap &\wto\nabla\ckap &&\text{in }\LSploc{4}{\bomega\times[0,\infty)},\label{eq:w-conv-nabla-cepkap}\\
\uepkap&\to\ukap &&\text{in }\LSploc{2}{\bomega\times[0,\infty)}\text{ and a.e. in }\Omega\times(0,\infty),\label{eq:conv-uepkap}\\
\uepkap&\wto\ukap &&\text{in }\LSploc{\frac{10}3}{\bomega\times[0,\infty)},\label{eq:w-conv-uepkap}\\
\nabla\uepkap&\wto\nabla\ukap &&\text{in }\LSploc{2}{\bomega\times[0,\infty)}\label{eq:w-conv-nabla-uepkap},
\end{alignat}
as $\epsi=\epsi_j\searrow0$. Moreover, the triple $(\nkap,\ckap,\ukap)$ is a global weak solution of \eqref{CTNScons},\eqref{BC} and \eqref{IC} in the sense of Definition \ref{def:sol}.
\end{lemma}

\begin{bew}
Combining the bounds of Lemmas \ref{lem:timespace-bounds} and \ref{lem:timereg} with an Aubin--Lions type lemma (\cite[Corollary 8.4]{Sim87}) we obtain that for any $\kappa\in[-1,1]$
\begin{align*}
\big\{\nepkap\big\}_{\epsi\in(0,1)}\quad\text{is relatively compact in }\LSploc{\frac{5}{4}}{\bomega\times[0,\infty)}
\end{align*}
and that hence there is some sequence $(\epsi_j)_{j\in\N}$ with $\epsi_j\searrow0$ as $j\to\infty$ such that $\nepkap\to\nkap$ in $\LSploc{\frac54}{\bomega\times[0,\infty)}$ and a.e. in $\Omega\times(0,\infty)$. According to the spatio-temporal bounds in Lemma \ref{lem:timespace-bounds} we can furthermore conclude \eqref{eq:w-conv-nabla-nepkap} and \eqref{eq:w-conv-nepkap} along a subsequence (which we still denote by $\epsi_j$). Moreover, also by Lemma \ref{lem:timespace-bounds}, $\{(\nepkap)^p\}_{\epsi\in(0,1)}$ is equi-integrable for any $p<\frac{5}{3}$ and therefore the a.e. convergence of $\nepkap$ together with Vitali's convergence theorem entail the strong convergence in \eqref{eq:conv-nepkap}. In a similar fashion we can make use of the bounds for $\cepkap$ in Lemmas \ref{lem:locex}, \ref{lem:bounds} and \ref{lem:timereg} to obtain \eqref{eq:conv-cepkap}--\eqref{eq:w-conv-nabla-cepkap} and the bounds for $\uepkap$ from Lemmas \ref{lem:bounds}, \ref{lem:timespace-bounds} and \ref{lem:timereg} to verify \eqref{eq:conv-uepkap}--\eqref{eq:w-conv-nabla-uepkap} upon extraction of another subsequence. That $(\nkap,\ckap,\ukap)$ solves \eqref{CTNScons} weakly in $\Omega\times(0,\infty)$ is then a straightforward consequence of the regularity and convergence properties we established just now, as these allow us to pass to the limit in all integrals making up the weak formulation of a solution, where we note that in particular \eqref{eq:conv-nepkap} and \eqref{eq:w-conv-uepkap} entail that for $\varphi\in C_0^\infty(\bomega\times[0,\infty))$ $\int_0^\infty\!\intomega\nepkap\uepkap\cdot\varphi\to\int_0^\infty\!\intomega\nkap\ukap\cdot\varphi$ and that \eqref{eq:conv-uepkap} and the dominated convergence theorem imply that $Y_\epsi\uepkap\to \ukap$ in $\LSploc{2}{\bomega\times[0,\infty)}$.
\end{bew}

\setcounter{equation}{0} 
\section{Eventual smallness of oxygen concentration with waiting times independent of \texorpdfstring{$\epsi$ and $\kappa$}{eps and kappa}}\label{sec4:decay}

The main objective of this section will be to establish several eventual smallness results for the chemical concentration, where, most importantly, the necessary waiting time of each estimate is independent of $\epsi\in(0,1)$ and $\kappa\in[-1,1]$. While it is known that these stabilizations occur in the setting with fixed $\kappa=1$ (\cite{win_chemonavstokesfinal_TransAm17}), the methods behind these results cannot be transferred directly if we want to maintain independence of the waiting time from the parameters $\epsi$ and $\kappa$. We start with two rather mild eventual smallness properties akin to \cite[Lemma 4.2]{win_chemonavstokesfinal_TransAm17}.

\begin{lemma}
\label{lem:decay-nabc-clnn}
For all $\delta>0$ there exists $T>0$ such that for each $\epsi\in(0,1)$ and $\kappa\in[-1,1]$ the solution $(\nepkap,\cepkap,\uepkap)$ of \eqref{approxprop} satisfies
\begin{align*}
\inf_{t\in[0,T]}\int_{t}^{t+1}\!\!\intomega\frac{1}{\epsi}\ln\big(1+\epsi\nepkap\big)\cepkap<\delta,
\end{align*}
as well as
\begin{align*}
\inf_{t\in[0,T]}\int_{t}^{t+1}\!\!\intomega\big|\nabla\cepkap\big|^2<\delta.
\end{align*}
\end{lemma}

\begin{bew}
Given $\delta>0$ we pick $T\in\N$ satisfying $\big(\|c_0\|_{\Lo[\infty]}+\|c_0\|_{\Lo[\infty]}^2\big)|\Omega|\delta^{-1}<T$. Then, utilizing the second and fourth equations of \eqref{approxprop} and the prescribed boundary conditions we find that for all $\epsi\in(0,1)$ and all $\kappa\in[-1,1]$ the equality
\begin{align*}
\frac{\intd}{\intd t}\intomega c_{\epsi t}^{(\kappa)}=\intomega\Delta\cepkap-\intomega\frac{1}{\epsi}\ln\big(1+\epsi\nepkap\big)\cepkap-\intomega\uepkap\cdot\nabla\cepkap=-\intomega\frac{1}{\epsi}\ln\big(1+\epsi \nepkap\big)\cepkap
\end{align*}
is valid on $(0,\infty)$. Integration over $(0,T)$ thus shows
\begin{align}\label{eq:decay-nabc-ineq1}
\intoTomega\frac{1}{\epsi}\ln\big(1+\epsi\nepkap\big)\cepkap\leq\intomega c_0\quad\text{for all }\epsi\in(0,1),\  \kappa\in[-1,1],
\end{align}
due to $\cepkap$ being nonnegative. Similarly, considering $\frac12\frac{\intd}{\intd t}\intomega(\cepkap)^2$ and making use of the fact that $\tfrac1{\epsi}\ln(1+\epsi s)\geq 0$ for all $\epsi\in(0,1)$ and $s\geq0$, we find that
\begin{align}\label{eq:decay-nabc-ineq2}
\intoTomega\big|\nabla\cepkap\big|^2\leq\intomega c_0^2\quad\text{for all }\epsi\in(0,1),\ \kappa\in[-1,1].
\end{align}
From \eqref{eq:decay-nabc-ineq1}, \eqref{eq:decay-nabc-ineq2} and Lemma \ref{lem:locex} we first obtain that for all $\epsi\in(0,1)$ and $\kappa\in[-1,1]$ we have
\begin{align*}
&\sum_{t=0}^{T-1}\Big(\int_t^{t+1}\!\intomega\frac{1}{\epsi}\ln\big(1+\epsi\nepkap\big)\cepkap+\int_t^{t+1}\!\intomega\big|\nabla\cepkap\big|^2\Big)\\
=\ &
\intoTomega\frac{1}{\epsi}\ln\big(1+\epsi\nepkap\big)\cepkap+\intoTomega\big|\nabla\cepkap\big|^2\leq \big(\|c_0\|_{\Lo[\infty]}+\|c_0\|_{\Lo[\infty]}^2\big)|\Omega|=:M
\end{align*}
and infer from this that for all $\epsi\in(0,1)$ and $\kappa\in[-1,1]$ there exists some $t_0\in[0,T]$ satisfying
\begin{align*}
\int_{t_0}^{t_0+1}\!\intomega\frac{1}{\epsi}\ln\big(1+\epsi\nepkap\big)\cepkap+\int_{t_0}^{t_0+1}\!\intomega\big|\nabla\cepkap\big|^2\leq \frac{M}{T}<\delta.
\end{align*}
In conclusion, for all $\delta>0$ one can find $T>0$ such that for all $\epsi\in(0,1)$ and $\kappa\in[-1,1]$ 
\begin{align*}
\inf_{t\in[0,T]}\Big(\int_{t}^{t+1}\!\intomega\frac{1}{\epsi}\ln\big(1+\epsi\nepkap\big)\cepkap+\int_{t}^{t+1}\!\intomega\big|\nabla\cepkap\big|^2\Big)<\delta,
\end{align*}
which clearly implies the assertion of the lemma.
\end{bew}

Making use of the uniform bounds from the previous sections and the lemma above we can also derive an additional eventual smallness property, which resembles the result of \cite[Lemma 4.4]{win_chemonavstokesfinal_TransAm17}.

\begin{lemma}
\label{lem:decay-c-L1}
For all $\delta>0$ there exists $T>0$ such that for any $\epsi\in(0,1)$ and $\kappa\in[-1,1]$ the solution $(\nepkap,\cepkap,\uepkap)$ of \eqref{approxprop} satisfies
\begin{align*}
\inf_{t\in[0,T]}\int_{t}^{t+1}\!\!\intomega\cepkap<\delta.
\end{align*}
\end{lemma}

\begin{bew}
As previously employed in the proof of Lemma \ref{lem:bounds}, we first note that the \GNI\ provides $C_1>0$ such that
\begin{align}\label{eq:decayL1-GNI}
\|\varphi\|_{\Lo[\frac{12}{5}]}^4\leq C_1\|\nabla\varphi\|_{\Lo[2]}\|\varphi\|_{\Lo[2]}^3+C_1\|\varphi\|_{\Lo[2]}^4
\end{align}
holds for all $\varphi\in\W[1,2]$. Moreover, the embedding $\W[1,2]\hookrightarrow\Lo[6]$ as well as the Poincaré inequality entail the existence of $C_2>0$ satisfying
\begin{align}\label{eq:decayL1-EmbPoin}
\|\varphi-\overline{\varphi}\|_{\Lo[6]}\leq C_2\|\nabla\varphi\|_{\Lo[2]}\quad\text{for all }\varphi\in\W[1,2],
\end{align}
where here and below we denote by $\overline{\varphi}:=\frac{1}{|\Omega|}\intomega\varphi$ the spatial average. Preparing later estimates we abbreviate $m:=\intomega n_0$ and set $C_3:=\frac{1}{2}\min\big\{|\Omega|,m\big\}$ and given any $\delta>0$ we then fix 
\begin{align}\label{eq:decayL1-delta0}
0<\delta_0<\min\left\{\frac{C_3\delta}{2|\Omega|},\frac{C_3^2\delta^2}{4|\Omega|^2C_1C_2^2m^\nfrac{3}{2}\big(K_1^\nfrac12+m^\nfrac12\big)}\right\},
\end{align}
where $K_1>0$ is the constant obtained in Lemma \ref{lem:bounds}. According to Lemma \ref{lem:decay-nabc-clnn}, one can find $T>0$ such that for any fixed $\epsi\in(0,1)$ and $\kappa\in[-1,1]$ there is some $t_0\in[0,T]$ satisfying
\begin{align}\label{eq:decayL1-estdelta0}
\int_{t_0}^{t_0+1}\!\!\intomega\frac{1}{\epsi}\ln\big(1+\epsi\nepkap\big)\cepkap<\delta_0\quad\text{and}\quad\int_{t_0}^{t_0+1}\!\!\intomega\big|\nabla\cepkap\big|^2<\delta_0.
\end{align}
To show that in fact this $T>0$ already fulfills the asserted property we continue by recalling that $\frac{1}{\epsi}\ln\big(1+\epsi s\big)\geq\frac{1}{2}\min\{s,1\}$ for all $\epsi\in(0,1)$ and $s\geq0$ and then estimate
\begin{align}\label{eq:decayL1-estL1}
\int_{t_0}^{t_0+1}\!\intomega&\,\frac{1}{\epsi}\ln\big(1+\epsi\nepkap\big)\cepkap-\int_{t_0}^{t_0+1}\!\intomega\frac{1}{\epsi}\ln\big(1+\epsi\nepkap\big)\big(\cepkap-\overline{c}_\epsi^{(\kappa)}\big)\nonumber\\
=\ &\int_{t_0}^{t_0+1}\overline{c}_\epsi^{(\kappa)}\intomega\frac{1}{\epsi}\ln\big(1+\epsi\nepkap\big)
\geq \int_{t_0}^{t_0+1}\overline{c}_\epsi^{(\kappa)}\intomega\frac{1}{2}\min\big\{\nepkap,1\big\}=\frac{C_3}{|\Omega|}\int_{t_0}^{t_0+1}\!\intomega\cepkap.
\end{align}
Making use of the Hölder inequality twice and drawing on \eqref{eq:decayL1-EmbPoin}, as well as the fact that $\frac1\epsi\ln(1+\epsi s)\leq s$ for all $\epsi\in(0,1)$ and $s\geq0$, we see that
\begin{align*}
-&\int_{t_0}^{t_0+1}\!\!\intomega\frac{1}{\epsi}\ln\big(1+\epsi\nepkap\big)\big(\cepkap-\overline{c}_\epsi^{(\kappa)}\big)\\&\quad\leq\left(\int_{t_0}^{t_0+1}\!\big\|\cepkap-\overline{c}_\epsi^{(\kappa)}\big\|^2_{\Lo[6]}\right)^\nfrac{1}{2}\left(\int_{t_0}^{t_0+1}\!\big\|\tfrac{1}{\epsi}\ln\big(1+\epsi\nepkap\big)\big\|_{\Lo[\frac65]}^2\right)^\nfrac{1}{2}\\&\qquad\leq C_2\left(\int_{t_0}^{t_0+1}\!\big\|\nabla\cepkap\big\|_{\Lo[2]}^2\right)^\nfrac{1}{2}\left(\int_{t_0}^{t_0+1}\!\big\|\nepkap\big\|_{\Lo[\frac{6}{5}]}^2\right)^\nfrac{1}{2}.
\end{align*}
Plugging this into \eqref{eq:decayL1-estL1} and combining with \eqref{eq:decayL1-estdelta0} therefore implies that
\begin{align*}
\int_{t_0}^{t_0+1}\!\intomega \cepkap\leq \frac{|\Omega|\delta_0}{C_3}+\frac{|\Omega|C_2\delta_0^\nfrac{1}{2}}{C_3}\left(\int_{t_0}^{t_0+1}\!\big\|\nepkap\big\|_{\Lo[\frac{6}{5}]}^2\right)^\nfrac{1}{2}.
\end{align*}
To further estimate the remaining term, we make use of \eqref{eq:decayL1-GNI}, the \CSI\ and Lemma \ref{lem:bounds} to find that
\begin{align*}
\ &\int_{t_0}^{t_0+1}\!\big\|\nepkap\big\|_{\Lo[\frac{6}{5}]}^2=\int_{t_0}^{t_0+1}\!\big\|{\nepkap}^\frac{1}{2}\big\|_{\Lo[\frac{12}{5}]}^4\leq C_1m^\nfrac{3}{2}\int_{t_0}^{t_0+1}\!\big\|\nabla{\nepkap}^\frac{1}{2}\big\|_{\Lo[2]}+C_1 m^2\\\leq\ &C_1m^\nfrac{3}{2}\left(\int_{t_0}^{t_0+1}\!\big\|\nabla{\nepkap}^\frac{1}{2}\big\|_{\Lo[2]}^2\right)^\nfrac{1}{2}+C_1m^2\leq C_1m^\nfrac{3}{2}\big(K_1^\nfrac{1}{2}+m^\nfrac{1}{2}\big),
\end{align*}
with $K_1>0$ provided by Lemma \ref{lem:bounds}. This, in light of \eqref{eq:decayL1-delta0}, establishes that
\begin{align*}
\int_{t_0}^{t_0+1}\!\intomega \cepkap\leq \frac{|\Omega|\delta_0}{C_3}+\frac{|\Omega|C_1^\nfrac{1}{2}C_2m^\nfrac{3}{4}\big(K_1^\nfrac{1}{2}+m^{\nfrac{1}{2}}\big)^\nfrac{1}{2}}{C_3}\delta_0^\nfrac{1}{2}<\delta,
\end{align*}
and thereby completes the proof. 
\end{bew}

Finally, augmenting the arguments of \cite[Lemma 3.4]{Lan-Longterm_M3AS16} to cover our setting, we obtain the eventual smallness of the $\Lo[\infty]$-norm of the oxygen concentration with waiting time uniform in $\epsi$ and $\kappa$.

\begin{lemma}
\label{lem:decay-c-Linfty}
For all $\delta>0$ there exists $T>0$ such that for each $\epsi\in(0,1)$, $\kappa\in[-1,1]$ and all $t>T$ the solution $(\nepkap,\cepkap,\uepkap)$ of \eqref{approxprop} satisfies
\begin{align*}
\big\|\cepkap(\cdot,t)\big\|_{\Lo[\infty]}<\delta.
\end{align*}
\end{lemma}

\begin{bew}
Similar to before we first note that by the \GNI\ we can find $C_1>0$ such that
\begin{align}\label{eq:decay-c-Linfty-GNI}
\|\varphi\|_{\Lo[\infty]}\leq C_1\|\nabla\varphi\|_{\Lo[4]}^\frac{12}{13}\|\varphi\|_{\Lo[1]}^\frac{1}{13}+C_1\|\varphi\|_{\Lo[1]}\quad\text{for all }\varphi\in\W[1,4].
\end{align}
Moreover, according to Lemma \ref{lem:bounds} there is $K_1>0$ such that
\begin{align}\label{eq:decay-c-Linfty-est1}
\int_t^{t+1}\!\!\intomega\big|\nabla\cepkap\big|^4\leq K_1
\end{align}
is valid for all $t>0$, $\epsi\in(0,1)$ and $\kappa\in[-1,1]$. Now, given $\delta>0$ we fix $0<\delta_0<\min\big\{\frac{\delta}{2C_1},\frac{\delta^{13}}{2^{13}C_1^{13}K_1^3}\big\}$ and note that in light of Lemma \ref{lem:decay-c-L1} we thus find $T_0>0$ such that for any fixed $\epsi\in(0,1)$ and $\kappa\in[-1,1]$ there is $t_0\in[0,T_0]$ satisfying
\begin{align}\label{eq:decay-c-Linfty-est2}
\int_{t_0}^{t_0+1}\!\!\intomega\cepkap<\delta_0.
\end{align}
From a combination of \eqref{eq:decay-c-Linfty-GNI} with two applications of Hölder's inequality, \eqref{eq:decay-c-Linfty-est1} and \eqref{eq:decay-c-Linfty-est2} we can directly conclude
\begin{align*}
\int_{t_0}^{t_0+1}\!\big\|\cepkap\big\|_{\Lo[\infty]}&\leq C_1\Big(\int_{t_0}^{t_0+1}\!\big\|\nabla\cepkap\big\|_{\Lo[4]}^4\Big)^\frac{3}{13}\Big(\int_{t_0}^{t_0+1}\!\big\|\cepkap\big\|_{\Lo[1]}\Big)^\frac{1}{13}+C_1\int_{t_0}^{t_0+1}\!\big\|\cepkap\big\|_{\Lo[1]}\\&\leq C_1K_1^\frac{3}{13}\delta_0^\frac{1}{13}+C_1\delta_0,
\end{align*}
which, by choice of $\delta_0$ implies
\begin{align*}
\int_{t_0}^{t_0+1}\!\big\|\cepkap\big\|_{\Lo[\infty]}<\delta.
\end{align*}
This entails that for all $\delta>0$ there exists $T_0>0$ such that for all $\epsi\in(0,1)$ and $\kappa\in[-1,1]$ one can find $t_0\in[0,T_0]$ such that
\begin{align*}
\inf_{t\in[t_0,t_0+1]}\big\|\cepkap(\cdot,t)\big\|_{\Lo[\infty]}\leq\int_{t_0}^{t_0+1}\!\big\|\cepkap\big\|_{\Lo[\infty]}<\delta,
\end{align*}
which, by recalling that $t\mapsto\big\|\cepkap(\cdot,t)\big\|_{\Lo[\infty]}$ is nonincreasing, immediately implies the assertion of the lemma with $T\geq T_0+1$.
\end{bew}

\setcounter{equation}{0} 
\section{Eventual \texorpdfstring{$L^p$}{Lp}-regularity estimates independent of \texorpdfstring{$\epsi$ and $\kappa$}{eps and kappa} as consequence of small oxygen concentration}\label{sec5:regest}

The uniform waiting time for smallness of $\cepkap$ in $\Lo[\infty]$ will be the key ingredient in obtaining additional regularity estimates for $\nepkap$ and $\uepkap$. We start by deriving a differential inequality for $\nepkap$ valid for all times after the chemical concentration has decayed below some threshold number $\eta$ which, in a second step, together with Lemma \ref{lem:decay-c-Linfty} will then show that the $\Lo[p]$-norm of $\nepkap$ is nonincreasing beyond some waiting time. A functional of similar form to the one we use in Lemma \ref{lem:testing} to derive the differential inequality has previously been successfully employed in e.g. \cite[Lemma 5.1]{win-stab2d-ArchRatMechAna12} and \cite[Lemma 3.5]{Lan-Longterm_M3AS16}.

\begin{lemma}
\label{lem:testing}
Let $T>0$, $p>1$, $\theta>0$ and $\eta>0$, $\epsi\in(0,1)$ and $\kappa\in[-1,1]$. If the solution $(\nepkap,\cepkap,\uepkap)$ of \eqref{approxprop} satisfies $$\big\|\cepkap(\cdot,t)\big\|_{\Lo[\infty]}\leq \eta\quad\text{for all }t>T,$$ then
\begin{align}\label{eq:testing-eq}
\frac{\intd}{\intd t}\intomega\frac{{\nepkap}^p}{\big(2\eta-\cepkap\big)^\theta}\leq &-p(p-1)\intomega\frac{{\nepkap}^{p-2}}{\big(2\eta-\cepkap\big)^\theta}\big|\nabla\nepkap\big|^2\nonumber\\
&\quad+\intomega\bigg(\frac{p(p-1)}{\big(1+\epsi\nepkap\big)\big(2\eta-\cepkap\big)^\theta}-\frac{2p\theta}{\big(2\eta-\cepkap\big)^{\theta+1}}\bigg){\nepkap}^{p-1}\big(\nabla\nepkap\cdot\nabla\cepkap\big)\\
&\qquad-\intomega\bigg(\frac{\theta(\theta+1)}{\big(2\eta-\cepkap\big)^{\theta+2}}-\frac{p\theta}{\big(1+\epsi\nepkap\big)\big(2\eta-\cepkap\big)^{\theta+1}}\bigg){\nepkap}^p\big|\nabla\cepkap\big|^2\nonumber
\end{align}
on $(T,\infty)$.
\end{lemma}

\begin{bew}
First we note that $t\mapsto\intomega\frac{(\nepkap(\cdot,t))^p}{(2\eta-\cepkap(\cdot,t))^\theta}$ is well-defined on $(T,\infty)$ due to $\big\|\cepkap(\cdot,t)\big\|_{\Lo[\infty]}\leq\eta$ for all $t>T$ and then a straightforward computation, utilizing integration by parts, shows
\begin{align}\label{eq:testing-eq-proof1}
\frac{\intd}{\intd t}\intomega\frac{{\nepkap}^p}{\big(2\eta-\cepkap\big)^\theta}=\;&p\intomega\frac{{\nepkap}^{p-1}}{\big(2\eta-\cepkap\big)^\theta}\bigg(\Delta\nepkap-\nabla\cdot\Big(\frac{\nepkap}{1+\epsi\nepkap}\nabla\cepkap\Big)-\nabla\nepkap\cdot\uepkap\bigg)\nonumber\\
&\quad+\theta\intomega\frac{{\nepkap}^p}{\big(2\eta-\cepkap\big)^{\theta+1}}\bigg(\Delta\cepkap-\frac{1}{\epsi}\ln\big(1+\epsi\nepkap\big)\cepkap-\nabla\cepkap\cdot\uepkap\bigg)\nonumber\\
\leq&-p(p-1)\intomega\frac{{\nepkap}^{p-2}\big|\nabla\nepkap\big|^2}{\big(2\eta-\cepkap\big)^\theta}-p\theta\intomega\frac{{\nepkap}^{p-1}\big(\nabla\nepkap\cdot\nabla\cepkap\big)}{\big(2\eta-\cepkap\big)^{\theta+1}}\\
&+p(p-1)\intomega\frac{{\nepkap}^{p-1}\big(\nabla\nepkap\cdot\nabla\cepkap\big)}{\big(1+\epsi\nepkap\big)\big(2\eta-\cepkap\big)^\theta}+p\theta\intomega\frac{{\nepkap}^p\big|\nabla\cepkap\big|^2}{\big(1+\epsi\nepkap\big)\big(2\eta-\cepkap\big)^{\theta+1}}\nonumber\\
&-p\theta\intomega\frac{{\nepkap}^{p-1}\big(\nabla\nepkap\cdot\nabla\cepkap\big)}{\big(2\eta-\cepkap\big)^{\theta+1}}-\theta(\theta+1)\intomega\frac{{\nepkap}^p\big|\nabla\cepkap\big|^2}{\big(2\eta-\cepkap\big)^{\theta+2}}\nonumber\\
&-\intomega\frac{\nabla(\nepkap)^p}{\big(2\eta-\cepkap\big)^\theta}\cdot\uepkap-\intomega{\nepkap}^p\nabla\big(2\eta-\cepkap\big)^{-\theta}\cdot\uepkap\nonumber
\end{align}
for all $t>T$, where we also made use of the fact that $\frac{1}{\epsi}\ln(1+\epsi s)\geq0$ for $s\geq0$. Herein, we have
\begin{align*}
-\intomega\frac{\nabla(\nepkap)^p}{\big(2\eta-\cepkap\big)^\theta}\cdot\uepkap-\intomega{\nepkap}^p\nabla\big(2\eta-\cepkap\big)^{-\theta}\cdot\uepkap=-\intomega\nabla\Big(\frac{{\nepkap}^p}{\big(2\eta-\cepkap\big)^\theta}\Big)\cdot\uepkap=0,
\end{align*}
due to the imposed boundary conditions and $\uepkap$ being divergence-free. Therefore, rearranging the terms of \eqref{eq:testing-eq-proof1} appropriately, we can immediately conclude \eqref{eq:testing-eq}.
\end{bew}

Waiting long enough for $\cepkap$ to decay past a certain threshold now entails the following.
\begin{lemma}
\label{lem:n-lp-step-bounds}
For all $p>1$ there exist $K_2>0$ and $T>0$ such that for any $\epsi\in(0,1)$, $\kappa\in[-1,1]$ and every $[t_1,t_2)\subseteq[T,\infty)$ the solution $(\nepkap,\cepkap,\uepkap)$ of \eqref{approxprop} satisfies
\begin{align*}
\intomega{\nepkap}^p(\cdot,t_2)+\int_{t_1}^{t_2}\!\intomega{\nepkap}^{p-2}\big|\nabla\nepkap\big|^2\leq K_2\intomega{\nepkap}^p(\cdot,t_1).
\end{align*}
\end{lemma}

\begin{bew}
Given $p>1$ we first fix $\theta\in(0,p-1)$ and then pick some $\eta>0$ satisfying 
\begin{align}\label{eq:n-lp-step-proof-eta}
\eta<\min\left\{\frac{\theta+1}{2p},\sqrt{\frac{\theta(\theta+1-\frac{p}{p-1}\theta)}{p(p-1)}}\right\}.
\end{align}
For these choices of parameters, in light of Lemma \ref{lem:decay-c-Linfty}, we find some $T=T(p)>0$ such that for all $\epsi\in(0,1)$ and $\kappa\in[-1,1]$ we have
\begin{align*}
\big\|\cepkap(\cdot,t)\big\|_{\Lo[\infty]}\leq \eta\quad\text{for all }t\geq T.
\end{align*}
Hence, the requirements of Lemma \ref{lem:testing} are met and the inequality \eqref{eq:testing-eq} is valid on $(T,\infty)$. Moreover, by choice of $\eta<\frac{\theta+1}{2p}$ and nonnegativity of $\nepkap$ and $\cepkap$ we have
\begin{align*}
\frac{\big(2\eta-\cepkap\big)p}{\big(1+\epsi\nepkap\big)(\theta+1)}\leq\frac{2\eta p}{\theta+1}<1\quad\text{on }[T,\infty)
\end{align*}
and hence
\begin{align*}
\frac{p\theta}{\big(1+\epsi\nepkap\big)\big(2\eta-\cepkap\big)^{\theta+1}}<\frac{\theta(\theta+1)}{\big(2\eta-\cepkap\big)^{\theta+2}}\quad\text{for all }t\geq T.
\end{align*}
Therefore, we can cancel out the term containing $\big|\nabla\cepkap\big|^2$ in \eqref{eq:testing-eq}. In fact, an employment of Young's inequality in \eqref{eq:testing-eq} shows that for all $\epsi\in(0,1)$ and $\kappa\in[-1,1]$ we have
\begin{align}\label{eq:n-lp-step-proof-diffineq}
\frac{\intd}{\intd t}\intomega\frac{{\nepkap}^p}{\big(2\eta-\cepkap\big)^\theta}\leq -\intomega\bigg(\frac{p(p-1)}{\big(2\eta-\cepkap\big)^\theta}-\frac{1}{4}H\big(\nepkap,\cepkap\big)\bigg){\nepkap}^{p-2}\big|\nabla\nepkap\big|^2
\end{align}
for all $t\geq T$, with
\begin{align*}
H(\sigma,\xi):=\frac{\Big(\frac{p(p-1)}{(1+\epsi\sigma)(2\eta-\xi)^\theta}-\frac{2p\theta}{(2\eta-\xi)^{\theta+1}}\Big)^2}{\frac{\theta(\theta+1)}{(2\eta-\xi)^{\theta+2}}-\frac{p\theta}{(1+\epsi\sigma)(2\eta-\xi)^{\theta+1}}},\quad\text{for }\sigma\geq0\text{ and }\xi\in[0,2\eta).
\end{align*}
To verify that in fact $\frac{p(p-1)}{(2\eta-\xi)^\theta}-\frac{1}{4}H(\sigma,\xi)\geq0$ for $\sigma\geq0$ and $\xi\in[0,2\eta)$, we first write
\begin{align*}
\frac{H(\sigma,\xi)(2\eta-\xi)^\theta}{4p(p-1)}=\frac{\frac{p(p-1)(2\eta-\xi)^2}{(1+\epsi\sigma)^2}-\frac{4p\theta(2\eta-\xi)}{1+\epsi\sigma}+\frac{4p\theta^2}{p-1}}{4\theta(\theta+1)-\frac{4p\theta(2\eta-\xi)}{1+\epsi\sigma}}=:\frac{H_1(\sigma,\xi)}{H_2(\sigma,\xi)}
\end{align*}
and note that by the nonnegativity of $\sigma$ and $\xi$ and latter part of \eqref{eq:n-lp-step-proof-eta} we have
\begin{align*}
H_1(\sigma,\xi)-H_2(\sigma,\xi)\leq p(p-1)4\eta^2+\frac{4p}{p-1}\theta^2-4\theta(\theta+1)<0.
\end{align*}
Since, due to \eqref{eq:n-lp-step-proof-eta}, we have $H_2(\sigma,\xi)\geq4\theta(\theta+1)-8p\theta\eta>0$ for $\sigma\geq0$ and $\xi\in[0,2\eta)$, this implies
\begin{align*}
\frac{H\big(\nepkap,\cepkap\big)\big(2\eta-\xi\big)^\theta}{4p(p-1)}\leq 1+\frac{p(p-1)4\eta^2+\frac{4p}{p-1}\theta^2-4\theta(\theta+1)}{4\theta(\theta+1)-8p\theta\eta}\quad\text{for all }\sigma\geq0,\ \xi\in[0,2\eta),
\end{align*}
from which we infer that
\begin{align*}
\frac{p(p-1)}{(2\eta-\xi)^\theta}-\frac{1}{4}H(\sigma,\xi)\geq C_3\frac{p(p-1)}{(2\eta-\xi)^\theta}>0\quad\text{for all }\sigma\geq0,\ \xi\in[0,2\eta),
\end{align*}
with $C_3:=-\frac{p(p-1)\eta^2+\frac{p}{p-1}\theta^2-\theta(\theta+1)}{\theta(\theta+1)-2p\eta\theta}>0$. Hence, we conclude from \eqref{eq:n-lp-step-proof-diffineq} that
\begin{align*}
\frac{\intd}{\intd t}\intomega\frac{{\nepkap}^p}{\big(2\eta-\cepkap\big)^\theta}+p(p-1)C_3\intomega\frac{{\nepkap}^{p-2}}{\big(2\eta-\cepkap\big)^\theta}\big|\nabla\nepkap\big|^2\leq 0\quad\text{for all }t\geq T,
\end{align*}
which for any $[t_1,t_2)\subseteq[T,\infty)$, upon integration with respect to time, shows that
\begin{align*}
\intomega\frac{{\nepkap}^p(\cdot,t_2)}{\big(2\eta-\cepkap(\cdot,t_2)\big)^\theta}+\int_{t_1}^{t_2}\!\intomega\frac{{\nepkap}^{p-2}}{\big(2\eta-\cepkap\big)^\theta}\big|\nabla\nepkap\big|^2\leq \frac{1}{\min\{1,p(p-1)C_3\}}\intomega\frac{{\nepkap}^p(\cdot,t_1)}{\big(2\eta-\cepkap(\cdot,t_1)\big)^\theta},
\end{align*}
completing the proof, after taking into account that $\eta^\theta\leq(2\eta-\cepkap)^\theta\leq (2\eta)^\theta$ on $\Omega\times[T,\infty)$.
\end{bew}

Making use of an inductive argument, as exercised in \cite[Lemma 6.3]{win_chemonavstokesfinal_TransAm17}, we can get rid of the time dependence present in the right-hand side of the inequality provided by Lemma \ref{lem:n-lp-step-bounds}. entailing the eventual uniform $L^p$-regularity of $\nepkap$ required for further analysis.

\begin{lemma}
\label{lem:n-lp-bounds}
For all $p>1$ there exist $T>0$ and $K_3=K_3(p)>0$ such that for each $\epsi\in(0,1)$, $\kappa\in[-1,1]$ and all $t>T$ the solution $(\nepkap,\cepkap,\uepkap)$ of \eqref{approxprop} satisfies
\begin{align*}
\intomega{\nepkap}^p(\cdot,t)\leq K_3\quad\text{and}\quad\int_{T}^\infty\!\intomega{\nepkap}^{p-2}\big|\nabla\nepkap\big|^2\leq K_3.
\end{align*}
\end{lemma}

\begin{bew}
Preparing an inductive argument we first assume that there exist $p_0>1$, $C_0>0$ and $T_0\geq0$ such that for all $\epsi\in(0,1)$, $\kappa\in[-1,1]$ and $t>T_0$ we have
\begin{align}\label{eq:n-lp-bounds-proof-eq1}
\int_{t}^{t+1}\big\|\nepkap\big\|_{\Lo[p_0]}\leq C_0.
\end{align}
In light of Lemma \ref{lem:n-lp-step-bounds} we find for each $q\in(1,p_0]$ corresponding $T_1=T_1(q)>0$ and $K_2=K_2(q)>0$ with the property that for all $\epsi\in(0,1)$, $\kappa\in[-1,1]$ and $[t_1,t)\subseteq[T_1,\infty)$ the inequality
\begin{align}\label{eq:n-lp-bounds-proof-eq2}
\intomega{\nepkap}^{q}(\cdot,t)+\int_{t_1}^{t}\!\intomega{\nepkap}^{q-2}\big|\nabla\nepkap\big|^2\leq K_2\intomega{\nepkap}^{q}(\cdot,t_1)
\end{align}
is valid. Letting $\bar{T}:=\max\{T_0,T_1\}$ we see that in view of \eqref{eq:n-lp-bounds-proof-eq1} there exists $C_1>0$ such that for any $\epsi\in(0,1)$ and $\kappa\in[-1,1]$ we can find $t_*\in[\bar{T},\bar{T}+1]$ such that $\|\nepkap(\cdot,t_*)\|_{\Lo[q]}\leq C_1$. Plugging this into \eqref{eq:n-lp-bounds-proof-eq2} with $t_1=t_*$ we obtain for all $t>\bar{T}+1$ and any $\epsi\in(0,1)$ and $\kappa\in[-1,1]$ that
\begin{align}\label{eq:n-lp-bounds-proof-eq3}
&\intomega{\nepkap}^{q}(\cdot,t)+\int_{\bar{T}+1}^{t}\!\intomega{\nepkap}^{p_0-2}\big|\nabla\nepkap\big|^2\nonumber\\\leq\ &\intomega{\nepkap}^{q}(\cdot,t)+\int_{t_*}^{t}\!\intomega{\nepkap}^{q-2}\big|\nabla\nepkap\big|^2\leq K_2\intomega{\nepkap}^{q}(\cdot,t_*)\leq K_2C_1,
\end{align}
proving that under the assumption \eqref{eq:n-lp-bounds-proof-eq1} the asserted bounds are valid for $p\in(1,p_0]$. Moreover, due to $\W[1,2]\hookrightarrow\Lo[6]$ and Hölder's inequality there is some $C_2=C_2(p_0)>0$ such that for all $t>\bar{T}+1$
\begin{align*}
\Big(\int_t^{t+1}\big\|\nepkap\big\|_{\Lo[3p_0]}\Big)^{p_0}&\leq \int_t^{t+1}\big\|{\nepkap}^{\frac{p_0}{2}}\big\|_{\Lo[6]}^2\leq C_2\int_t^{t+1}\Big(\big\|\nabla{\nepkap}^\frac{p_0}{2}\big\|_{\Lo[2]}^2+\big\|{\nepkap}^{\frac{p_0}{2}}\big\|_{\Lo[\frac{2}{p_0}]}^2\Big)\\
&\leq \frac{C_2 p_0^2}{4}\int_t^{t+1}\!\intomega{\nepkap}^{p_0-2}\big|\nabla\nepkap\big|+C_2 m^{p_0}\leq\frac{C_1C_2K_2p_0^2}{4}+C_2m^{p_0},
\end{align*}
where we also made use of $\intomega\nepkap(\cdot,t)=\intomega n_0=:m$ for all $t>0$ and \eqref{eq:n-lp-bounds-proof-eq3}. Drawing on these calculations the step from $p_0$ to $3p_0$ is  possible and we only have to ensure that indeed the assumption \eqref{eq:n-lp-bounds-proof-eq1} is fulfilled for some $p_0>1$. Now, in a similar fashion the embedding $\W[1,2]\hookrightarrow\Lo[6]$ and Lemma \ref{lem:bounds} provide $C_3>0$ and $K_1>0$ such that for all $\epsi\in(0,1)$, $\kappa\in[-1,1]$ and $t>0$ we have
\begin{align*}
\int_t^{t+1}\big\|\nepkap\big\|_{\Lo[3]}&\leq C_3\int_t^{t+1}\Big(\big\|\nabla{\nepkap}^\frac{1}{2}\|_{\Lo[2]}^2+\big\|{\nepkap}^\frac{1}{2}\big\|_{\Lo[2]}^2\Big)\\&\leq\frac{C_3}{4}\int_{t}^{t+1}\!\intomega\frac{\big|\nabla\nepkap\big|^2}{\nepkap}+C_3\int_t^{t+1}\!\intomega\nepkap\leq \frac{C_3K_1}{4}+C_3 m,
\end{align*}
which shows that \eqref{eq:n-lp-bounds-proof-eq3} is valid for $p_0=3$ and thereby concludes the proof.
\end{bew}

An immediate consequence is the eventual boundedness of the forcing term in the third equation of \eqref{approxprop}, from which we extract new regularity information on the gradient of $\uepkap$.

\begin{lemma}
\label{lem:ev-l2-bound-u}
There exist $T>0$ and $C>0$ such that for each $\epsi\in(0,1)$, $\kappa\in[-1,1]$ and all $t>T$ the solution $(\nepkap,\cepkap,\uepkap)$ of \eqref{approxprop} satisfies
\begin{align*}
\intomega\big|\nabla\uepkap(\cdot,t)\big|^2\leq C.
\end{align*}
\end{lemma}

\begin{bew}
Recalling that $\mathcal{P}$ denotes the Helmholtz projection from $\Lo[2]$ to $L^2_{\sigma}(\Omega)$ and $A:=-\mathcal{P}\Delta$ the Stokes operator in $\Lo[2]$ under homogeneous Dirichlet boundary conditions, we find that testing the projected third equation of \eqref{approxprop} by $A\uepkap$ implies in view of Young's inequality that
\begin{align}\label{eq:ev-l2-bound-u-eq1}
\frac{1}{2}\frac{\intd}{\intd t}\intomega\big|A^{\frac{1}{2}}\uepkap\big|^2+\intomega\big|A\uepkap\big|^2\leq\intomega\big|A\uepkap\big|^2+\frac{|\kappa|}{2}\intomega\Big|\big(Y_\epsi\uepkap\cdot\nabla\big)\uepkap\Big|^2+\frac{1}{2}\intomega\big|\nepkap\nabla\phi\big|^2
\end{align}
is valid for all $t>0$, where we also made use of the facts that $\big\|A^{\frac{1}{2}}\uepkap\big\|_{\Lo[2]}=\big\|\nabla\uepkap\big\|_{\Lo[2]}$ and $\|\mathcal{P}\varphi\|_{\Lo[2]}\leq\|\varphi\|_{\Lo[2]}$ for all $\varphi\in \Lo[2]$. Moreover, since $D(1+\epsi A)=\W[2,2]\cap W_{0,\sigma}^{1,2}\hookrightarrow\Lo[\infty]$ we see that $Y_\epsi\uepkap$ there exists $C_1>0$ such that for any $\epsi\in(0,1)$, $\kappa\in[-1,1]$ and all $t>0$ we have
\begin{align*}
\big\|Y_\epsi\uepkap(\cdot,t)\big\|_{\Lo[\infty]}\leq C_1\big\|\uepkap(\cdot,t)\big\|_{\Lo[2]}\leq C_1\sqrt{\frac{K_1}{K_0}}=:C_2,
\end{align*}
where $K_0,K_1>0$ are the constants obtained in Lemma \ref{lem:bounds}. In particular, we obtain from a combination with \eqref{eq:ev-l2-bound-u-eq1} and the fact that $|\kappa|\leq 1$ that
\begin{align*}
\frac{\intd}{\intd t}\intomega\big|\nabla\uepkap\big|^2\leq C_2\intomega\big|\nabla\uepkap\big|^2+\|\nabla\phi\|_{\Lo[\infty]}^2\intomega\big|\nepkap\big|^2
\end{align*}
holds for all $t>0$. Letting $y(t):=\intomega|\nabla\uepkap(\cdot,t)|^2$ and $h(t):=\|\nabla\phi\|_{\Lo[\infty]}^2\intomega|\nepkap(\cdot,t)|^2$ we find that by Lemma \ref{lem:n-lp-bounds} there exist $T>0$ and $C_3>0$ such that for any $\epsi\in(0,1)$, $\kappa\in[-1,1]$ and all $t\geq T$ we have $h(t)\leq C_3$, and hence
\begin{align}\label{eq:ev-l2-bound-u-eq2}
y'(t)\leq C_2 y(t)+C_3\quad\text{for all }t>T.
\end{align}
Recalling that with $K_1>0$ provided by Lemma \ref{lem:bounds} for any $\epsi\in(0,1)$, $\kappa\in[-1,1]$ and all $t>0$ we moreover have 
\begin{align*}
\int_t^{t+1}\!\intomega\big|\nabla\uepkap\big|^2\leq K_1,
\end{align*}
we infer that for any fixed $t>T+1$ and each $\epsi\in(0,1)$ and $\kappa\in[-1,1]$ there exists some $t_*\in(t-1,t)$ such that
\begin{align*}
\intomega\big|\nabla\uepkap(\cdot,t_*)\big|^2\leq K_1,
\end{align*}
which upon integrating the differential inequality \eqref{eq:ev-l2-bound-u-eq2} over $(t_*,t)$ shows that
\begin{align*}
y(t)\leq y(t_*)e^{C_2(t-t_*)}+\int_{t_*}^t C_3 e^{C_2(t-t_*)}\leq K_1e^{C_2}+C_3 e^{C_2},
\end{align*}
completing the proof.
\end{bew}

As last step in this section we also lift the regularity of the signal gradient for times beyond the waiting times from the previous lemmas.

\begin{lemma}
\label{lem:ev-l4-bound-nabc}
There exist $T>0$ and $C>0$ such that for each $\epsi\in(0,1)$, $\kappa\in[-1,1]$ and all $t>T$ the solution $(\nepkap,\cepkap,\uepkap)$ of \eqref{approxprop} satisfies
\begin{align*}
\intomega\big|\nabla\cepkap(\cdot,t)\big|^4\leq C.
\end{align*}
\end{lemma}

\begin{bew}
We work along similar lines as in the proof of Lemma \ref{lem:ev-l2-bound-u}, by first establishing a differential inequality for the quantity $\intomega\big|\nabla\cepkap(\cdot,t)\big|^4$. A standard testing procedure utilizing the pointwise identity $\nabla\cepkap\cdot\nabla\Delta\cepkap=\frac{1}{2}\Delta\big|\nabla\cepkap\big|-\big|D^2\cepkap\big|^2$, the fact that $\nabla\cdot\uepkap=0$ on $\Omega\times(0,\infty)$ and the upper estimate of \eqref{eq:approxest} shows that
\begin{align}\label{eq:ev-l4-bound-nabc-eq1}
\frac{1}{4}\frac{\intd}{\intd t}\intomega\big|\nabla\cepkap\big|^4&\leq\frac{1}{2}\intomega\big|\nabla\cepkap\big|^2\Delta\big|\nabla\cepkap\big|^2-\intomega\big|\nabla\cepkap\big|^2\big|D^2\cepkap\big|^2-\intomega\big|\nabla\cepkap\big|^2 \nabla\cepkap\cdot\big(\nabla\uepkap\cdot\nabla\cepkap\big)\nonumber\\
&\quad+\intomega\big|\nabla\cepkap\big|^2\Delta\cepkap\nepkap\cepkap+2\intomega\nabla\cepkap\cdot\big(D^2\cepkap\cdot\nabla\cepkap\big)\nepkap\cepkap
\end{align}
holds for all $t>0$. To further estimate the first term on the right, we draw on arguments employed in \cite[Proposition 3.2]{ISY-quasilin-pp-JDE14}. Let us recall that there exists $C_1>0$ such that for any $\varphi\in\CSp{2}{\bomega}$ we have $\frac{\partial|\nabla\varphi|^2}{\partial\nu}\leq C_1|\nabla\varphi|^2$ on $\romega$ (cf. \cite[Lemma 4.2]{MS14}). Moreover, by utilizing the fact that for $r\in(0,\frac{1}{2})$ $\W[1,2]\hookrightarrow\hookrightarrow\W[r+\frac{1}{2},2]\hookrightarrow\Lo[1]$ (\cite{MR2944369}), Ehrling's lemma, as well as trace embddings (e.g. \cite[Theorem 4.24, Proposition 4.22]{HT08}), for every fixed $\eta>0$ we obtain $C_2>0$ such that $\|\psi\|_{L^2(\romega)}\leq\eta\|\nabla\psi\|_{\Lo[2]}+C_2\|\psi\|_{\Lo[1]}$ holds for any $\psi\in\W[1,2]$. Hence, drawing on Lemmas \ref{lem:locex} and \ref{lem:bounds} to estimate $\intomega|\nabla\cepkap|^2\leq 2K_1\|c_0\|_{\Lo[\infty]}=:C_3$, we find that for any $\epsi\in(0,1)$, $\kappa\in[-1,1]$ and all $t>0$
\begin{align*}
\frac{1}{2}\intomega\big|\nabla\cepkap\big|^2\Delta\big|\nabla\cepkap\big|^2&=-\frac{1}{2}\intomega\Big|\nabla\big|\nabla\cepkap\big|^2\Big|^2+\frac{1}{2}\intromega\big|\nabla\cepkap\big|^2\frac{\partial\big|\nabla\cepkap\big|^2}{\partial\nu}\\
&\leq -\frac{1}{2}\intomega\Big|\nabla\big|\nabla\cepkap\big|	^2\Big|^2+\frac{1}{2}\intomega\Big|\nabla\big|\nabla\cepkap\big|^2\Big|^2+C_2\Big(\intomega\big|\nabla\cepkap\big|^2\Big)^2
\leq C_2 C_3^2
\end{align*}
holds. Combining this with \eqref{eq:ev-l4-bound-nabc-eq1}, multiple employments of Young's inequality show that
\begin{align*}
\frac{1}{4}\frac{\intd}{\intd t}\intomega\big|\nabla\cepkap\big|^4
&\leq -\intomega\big|\nabla\cepkap\big|^2\big|D^2\cepkap\big|^2+\intomega\big|\nabla\cepkap\big|^4\big|\nabla\uepkap\big|+\frac{1}{12}\intomega\big|\nabla\cepkap\big|^2\big|\Delta\cepkap\big|^2\\&\quad+3\intomega\big|\nabla\cepkap\big|^2{\nepkap}^2{\cepkap}^2+\frac{1}{4}\intomega\big|\nabla\cepkap\big|^2\big|D^2\cepkap\big|^2+4\intomega\big|\nabla\cepkap\big|^2{\nepkap}^2{\cepkap}^2+C_2 C_3^2
\end{align*}
is valid for $t>0$. In light of the pointwise estimate $\big|\Delta\cepkap\big|^2\leq 3\big|D^2\cepkap\big|^2$ and Hölder's inequality this implies that
\begin{align}\label{eq:ev-l4-bound-nabc-eq2}
\frac{1}{4}\frac{\intd}{\intd t}\intomega\big|\nabla\cepkap\big|^4&+\frac{1}{2}\intomega\big|\nabla\cepkap\big|^2\big|D^2\cepkap\big|^2\\&\leq \Big(\intomega\big|\nabla\cepkap\big|^8\Big)^\frac{1}{2}\Big(\intomega\big|\nabla\uepkap\big|^2\Big)^\frac{1}{2}+7\Big(\intomega{\nepkap}^4{\cepkap}^4\Big)^\frac{1}{2}\Big(\intomega\big|\nabla\cepkap\big|^4\Big)^\frac{1}{2}+C_2 C_3^2\nonumber
\end{align}
for all $t>0$. Making use of the \GNI\ we obtain $C_4>0$ such that for each $\epsi\in(0,1)$, $\kappa\in[-1,1]$ and all $t>0$
\begin{align*}
\Big(\intomega\big|\nabla\cepkap\big|^8\Big)^\frac{1}{2}=\Big\|\big|\nabla\cepkap\big|^2\Big\|_{\Lo[4]}^2&\leq C_4\Big\|\nabla\big|\nabla\cepkap\big|^2\Big\|_{\Lo[2]}^{\frac{9}{5}}\Big\|\big|\nabla\cepkap\big|^2\Big\|_{\Lo[1]}^\frac{1}{5}+C_4\Big\|\big|\nabla\cepkap\big|^2\Big\|_{\Lo[1]}^2\\&\leq C_4 C_3^\frac{1}{5}\Big\|\nabla\big|\nabla\cepkap\big|^2\Big\|_{\Lo[2]}^{\frac{9}{5}}+C_4 C_3^2\\&\leq C_4 C_3^\frac{1}{5}\Big(4\intomega\big|\nabla\cepkap|^2 \big|D^2\cepkap\big|^2\Big)^\frac{9}{10}+C_4 C_3^2.
\end{align*}
Therefore, again by using Young's inequality, we infer from \eqref{eq:ev-l4-bound-nabc-eq2} that there is $C_5>0$ such that for each $\epsi\in(0,1)$, $\kappa\in[-1,1]$ and all $t>0$
\begin{align*}
\frac{1}{4}\frac{\intd}{\intd t}\intomega\big|\nabla\cepkap\big|^4\leq C_5\Big(\intomega\big|\nabla\uepkap\big|^2\Big)^5+C_4 C_3^2\Big(\intomega\big|\nabla\uepkap\big|^2\Big)^\frac{1}{2}+7\intomega{\nepkap}^4{\cepkap}^4+\frac{7}{4}\intomega\big|\nabla\cepkap\big|^4+C_4 C_3^2,
\end{align*}
which in light of Lemmas \ref{lem:locex}, \ref{lem:n-lp-bounds} and \ref{lem:ev-l2-bound-u} entails that there are $T>0$ and $C_6>0$ such that for each $\epsi\in(0,1)$, $\kappa\in[-1,1]$ and all $t>T$ the function $y(t):=\intomega\big|\nabla\cepkap(\cdot,t)\big|^4$ satisfies the differential inequality
\begin{align}\label{eq:ev-l4-bound-nabc-eq3}
y'(t)\leq 7y(t)+C_6.
\end{align}
Now, according to Lemmas \ref{lem:locex} and \ref{lem:bounds} for any $\epsi\in(0,1)$, $\kappa\in[-1,1]$ and all $t>0$ we can estimate
\begin{align*}
\int_t^{t+1}\!\intomega\big|\nabla\cepkap\big|^4\leq\|c_0\|_{\Lo[\infty]}^3\int_t^{t+1}\!\intomega\frac{\big|\nabla\cepkap\big|^4}{{\cepkap}^3}\leq \|c_0\|_{\Lo[\infty]}^3 K_1,
\end{align*} 
and hence for any fixed $t>T+1$ and each $\epsi\in(0,1)$ and $\kappa\in[-1,1]$ we find $t_*\in(t-1,t)$ such that
\begin{align*}
\intomega\big|\nabla\cepkap(\cdot,t_*)\big|^4\leq \|c_0\|_{\Lo[\infty]}^3 K_1=:C_7,
\end{align*}
which upon integrating \eqref{eq:ev-l4-bound-nabc-eq3} over $(t_*,t)$ entails that
\begin{align*}
y(t)\leq C_7e^7+C_6e^7
\end{align*}
as desired.
\end{bew}

\setcounter{equation}{0} 
\section{Uniform eventual stabilization of \texorpdfstring{$\nepkap$ and $\uepkap$ in some $L^p$}{nepkap and uepkap in some Lp} spaces}\label{sec6:ev-dec-n-u}

Eventual decay of the signal component and uniform regularity estimates at hand, we can now turn towards obtaining eventual stabilization properties of the two remaining solution components. These will be an important cornerstone of the maximal Sobolev regularity type arguments we employ in Section \ref{sec7:evsmooth} to obtain uniform bounds in Hölder spaces. We start with an eventual smallness result for a mixed quantity of $\nepkap$ and $\nabla\cepkap$.

\begin{lemma}
\label{lem:decay-grad-c-L4}
For all $\delta>0$ there exists $T>0$ such that for any $\epsi\in(0,1)$ and $\kappa\in[-1,1]$ and all $t>T$ the solution $(\nepkap,\cepkap,\uepkap)$ of \eqref{approxprop} satisfies
\begin{align*}
\int_t^{t+1}\big\|\nepkap\nabla\cepkap\big\|^2_{\Lo[2]}<\delta.
\end{align*}
\end{lemma}

\begin{bew}
According to Lemma \ref{lem:bounds} there is $K_1>0$ such that for each $\epsi\in(0,1)$, $\kappa\in[-1,1]$ and all $t>0$ we have 
\begin{align*}
\int_t^{t+1}\!\intomega\frac{\big|\nabla\cepkap\big|^4}{{\cepkap}^3}\leq K_1.
\end{align*}
Similarly, drawing on Lemma \ref{lem:n-lp-bounds}, we find $K_3>0$ and $T_1>0$ such that for any $\epsi\in(0,1)$, $\kappa\in[-1,1]$ and all $t>T_1$ the estimate
\begin{align*}
\big\|\nepkap(\cdot,t)\big\|_{\Lo[4]}^4\leq K_3
\end{align*}
is valid. Now, given any $\delta>0$ we fix 
\begin{align*}
0<\delta_0<\min\left\{\frac{\delta}{2K_1},\sqrt{\frac{\delta}{2K_3}}\right\}
\end{align*}
and, according to Lemma \ref{lem:decay-c-Linfty}, obtain a corresponding $T_2>0$ such that for each $\epsi\in(0,1)$, $\kappa\in[-1,1]$ and all $t>T_2$
\begin{align*}
\big\|\cepkap(\cdot,t)\big\|_{\Lo[\infty]}<\delta_0
\end{align*}
is satisfied. Hence, by making use of the estimates above, as well as Hölder's and Young's inequalities, we achieve for any $\epsi\in(0,1)$, $\kappa\in[-1,1]$ and all $t>\max\{T_1,T_2\}$ that
\begin{align*}
\int_t^{t+1}\big\|\nepkap\nabla\cepkap\big\|_{\Lo[2]}^2&\leq\int_t^{t+1}\big\|\nepkap\big\|_{\Lo[4]}^2\big\|\nabla\cepkap\big\|_{\Lo[4]}^2\leq\int_t^{t+1}K_3^\frac{1}{2}\big\|\nabla\cepkap\big\|_{\Lo[4]}^2\\
&\leq\int_t^{t+1}K_3^\frac{1}{2}\big\|\cepkap\big\|_{\Lo[\infty]}^\frac{3}{2}\left(\intomega\frac{\big|\nabla\cepkap\big|^4}{{\cepkap}^3}\right)^\frac{1}{2}\\
&\leq\int_t^{t+1}K_3\big\|\cepkap\big\|_{\Lo[\infty]}^2+\int_t^{t+1}\big\|\cepkap\big\|_{\Lo[\infty]}\intomega\frac{\big|\nabla\cepkap\big|^4}{{\cepkap}^3}\\&
\leq K_3\delta_0^2+K_1\delta_0<\delta,
\end{align*}
completing the proof.
\end{bew}

In order to successfully extract a uniform stabilization for $\nepkap$ and $\uepkap$ in certain $L^p$ spaces we will require the following auxiliary lemma for ODEs, which we have taken from \cite[Lemma 4.3]{WangWinklerXiang-smallconvection-MathZ18}.

\begin{lemma}
\label{lem:ODI-helplemma}
Let $I$ be any set and $\lambda>0$, and for each $\iota\in I$ let $y_\iota\in\CSp{0}{[0,\infty)}\cap\CSp{1}{(0,\infty)}$ and $f_\iota\in\CSp{0}{(0,\infty)}$ be nonnegative and such that
\begin{align*}
y'_\iota(t)+\lambda y_\iota(t)\leq f_\iota(t)\quad\text{for all }t>0
\end{align*}
and
\begin{align*}
\sup_{\iota\in I} y_\iota(0)<\infty,\quad\text{as well as}\quad
\sup_{\iota\in I}\|f_\iota\|_{\LSp{\infty}{(0,\infty)}}<\infty, \quad\text{and}\quad
\sup_{\iota\in I}\int_t^{t+1} f_\iota(s)\intd s\to 0\quad\text{as }t\to\infty
\end{align*}
hold. Then
\begin{align*}
\sup_{\iota\in I}y_\iota(t)\to0\quad\text{as }t\to\infty.
\end{align*}
\end{lemma}

Tracking the time evolution of $y_\epsi^{(\kappa)}(t):=\displaystyle\intomega \big(\nepkap(\cdot,t)-\overline{n_0}\big)^2$ and shifting the time appropriately, we can make use of the statement above to attain a uniform eventual smallness of $y_\epsi^{(\kappa)}$.

\begin{lemma}
\label{lem:nepkap-decay}
For all $\delta>0$ there exists $T>0$ such that for any $\epsi\in(0,1)$, $\kappa\in[-1,1]$ and all $t>T$ he solution $(\nepkap,\cepkap,\uepkap)$ of \eqref{approxprop} satisfies
\begin{align*}
\big\|\nepkap(\cdot,t)-\overline{n_0}\big\|_{\Lo[2]}^2< \delta.
\end{align*}
Furthermore, for all $p\geq2$ and $\delta'>0$ there is $T'>0$ such that for each $\epsi\in(0,1)$, $\kappa\in[-1,1]$ and all $t>T'$ the solution satisfies
\begin{align*}
\int_t^{t+1}\big\|\nepkap-\overline{n_0}\big\|_{\Lo[p]}^p<\delta'.
\end{align*}
\end{lemma}

\begin{bew}
We start with the case $p=2$. Due to the Young and Poincaré inequalities we obtain $C_1>0$ such that
\begin{align*}
\frac{\intd}{\intd t}\intomega\big(\nepkap-\overline{n_0}\big)^2\leq-\intomega\big|\nabla\nepkap\big|^2+\intomega\big|\nepkap\nabla\cepkap\big|^2\leq-\frac{1}{C_1}\intomega\big(\nepkap-\overline{n_0}\big)^2+\big\|\nepkap\nabla\cepkap\big\|_{\Lo[2]}^2
\end{align*}
holds for all $t>0$, where we also made use of the fact that $\nabla\cdot\uepkap=0$ in $\Omega\times(0,\infty)$. Moreover, as 
\begin{align*}
\big\|\big(\nepkap\nabla\cepkap\big)(\cdot,t)\big\|_{\Lo[2]}^2\leq\big\|\nepkap(\cdot,t)\big\|_{\Lo[4]}^2\big\|\nabla\cepkap(\cdot,t)\big\|_{\Lo[4]}^2\quad\text{for all }t>0,
\end{align*}
in light of Lemmas \ref{lem:n-lp-bounds}, \ref{lem:ev-l4-bound-nabc} and \ref{lem:decay-grad-c-L4}, we find that there exists some $T_1>0$ such that $$y_\epsi^{(\kappa)}(t):=\intomega\big(\nepkap(t-T_1)-\overline{n_0}\big)^2$$ satisfies the conditions of Lemma \ref{lem:ODI-helplemma}, and hence we conclude that for all $\delta>0$ there exists $\bar{T}\geq T_1$ such that for any $\epsi\in(0,1)$, $\kappa\in[-1,1]$ and all $t>\bar{T}$ we have
\begin{align*}
\big\|\nepkap(\cdot,t)-\overline{n_0}\big\|_{\Lo[2]}^2<\delta,
\end{align*} 
which in particular also immediately implies the second claim for $p=2$. For $p>2$ we let $K_3:=K_3(2p)>0$ and $T_2>0$ be given by Lemma \ref{lem:n-lp-bounds} and then, in consideration of \eqref{IR}, obtain $C_2>0$ such that for each $\epsi\in(0,1)$, $\kappa\in[-1,1]$ and all $t>T_2$ we have
\begin{align*}
\big\|\nepkap(\cdot,t)-\overline{n_0}\big\|_{\Lo[2p]}^\frac{p(p-2)}{p-1}\leq (K_3+\|\overline{n_0}\|_{\Lo[2p]})^\frac{p(p-2)}{p-1}\leq C_2.
\end{align*} 
Therefore, by means of Hölder interpolation and Hölder's inequality we find that for any $\epsi\in(0,1)$, $\kappa\in[-1,1]$ and all $t>T_2$
\begin{align}\label{eq:nepkap-decay-eq1}
\int_t^{t+1}\big\|\nepkap-\overline{n_0}\big\|_{\Lo[p]}^p&\leq\int_t^{t+1}\big\|\nepkap-\overline{n_0}\big\|_{\Lo[2p]}^\frac{p(p-2)}{p-1}\big\|\nepkap-\overline{n_0}\big\|_{\Lo[2]}^\frac{p}{p-1}\nonumber\\
&\leq C_2\int_t^{t+1}\big\|\nepkap-\overline{n_0}\big\|_{\Lo[2]}^\frac{p}{p-1}\leq C_2\int_t^{t+1}\big\|\nepkap-\overline{n_0}\big\|_{\Lo[2]}^2
\end{align}
is valid, due to $\frac{p}{p-1}<2$. Finally, for given $\delta>0$ we let $0<\delta_0<\frac{\delta}{C_2}$ and then conclude the proof by making use of the first part of this Lemma to estimate the remaining term in \eqref{eq:nepkap-decay-eq1} by $\delta_0$ for $t>T_3$ large enough.
\end{bew}

The second conclusion we can draw from the ODE lemma \ref{lem:ODI-helplemma} concerns the gradient of the fluid velocity field and, by Sobolev embeddings, the fluid velocity itself.

\begin{lemma}
\label{lem:uepkap-decay}
For all $\delta>0$ there exists $T>0$ such that for each $\epsi\in(0,1)$, $\kappa\in[-1,1]$ and all $t>T$ the solution $(\nepkap,\cepkap,\uepkap)$ of \eqref{approxprop} satisfies
\begin{align}\label{eq:uepkap-nab-decay}
\int_t^{t+1}\!\intomega\big|\nabla\uepkap\big|^2<\delta.
\end{align}
Moreover, for all $p\in[1,6]$ and all $\delta'>0$ there exists $T'>0$ such that for each $\epsi\in(0,1)$, $\kappa\in[-1,1]$ and all $t>T'$
\begin{align}\label{eq:uepkap-decay}
\int_t^{t+1}\big\|\uepkap\big\|_{\Lo[p]}^2<\delta'
\end{align}
holds.
\end{lemma}

\begin{bew}
Making use of Lemma \ref{lem:diff-ineq-2} and the divergence-free property of $\uepkap$ we first find that for each $\epsi\in(0,1)$, $\kappa\in[-1,1]$ and all $t>0$ we have
\begin{align*}
\frac{1}{2}\frac{\intd}{\intd t}\intomega\big|\uepkap\big|^2+\intomega\big|\nabla\uepkap\big|^2=\intomega\big(\nepkap-\overline{n_0}\big)\nabla\phi\cdot\uepkap.
\end{align*}
Here, the Poincaré inequality provides $C_1>0$ such that for each $\epsi\in(0,1)$, $\kappa\in[-1,1]$ we have $\|\uepkap\|_{\Lo[2]}\leq C_1\|\nabla\uepkap\|_{\Lo[2]}$, which entails upon an application of Young's inequality that for any $\epsi\in(0,1)$, $\kappa\in[-1,1]$ and all $t>0$
\begin{align}\label{eq:uepkap-decay-eq1}
\frac12\frac{\intd}{\intd t}\intomega\big|\uepkap\big|^2+\frac12\intomega\big|\nabla\uepkap\big|^2\leq \frac{C_1\big\|\nabla\phi\big\|_{\Lo[\infty]}^2}{2}\intomega\big(\nepkap-\overline{n_0}\big)^2
\end{align}
is valid on $(0,\infty)$. Since the Poincaré inequality moreover implies that for any $\epsi\in(0,1)$, $\kappa\in[-1,1]$ and all $t>0$
\begin{align*}
\frac{\intd}{\intd t}\intomega\big|\uepkap\big|^2+\frac{1}{C_1}\intomega\big|\uepkap\big|^2\leq C_1\|\nabla\phi\|_{\Lo[\infty]}^2\intomega\big(\nepkap-\overline{n_0}\big)^2
\end{align*}
holds, we find that in light of Lemmas \ref{lem:nepkap-decay}, \ref{lem:bounds}, \ref{lem:n-lp-bounds} and \eqref{phidef}, there exists some $T_1>0$ such that the function
\begin{align*}
y_\epsi^{(\kappa)}(t):=\intomega\big|\uepkap(t-T_1)\big|^2
\end{align*}
satisfies the conditions of Lemma \ref{lem:ODI-helplemma}. Hence, we find that for all $\delta_0>0$ there is some $\bar{T}\geq T_1$ such that for any $\epsi\in(0,1)$, $\kappa\in[-1,1]$ and all $t>\bar{T}$
\begin{align*}
\big\|\uepkap(\cdot,t)\big\|_{\Lo[2]}^2<\delta_0
\end{align*}
holds. Now, by making use of the first part of the proof and Lemma \ref{lem:nepkap-decay}, given any $\delta>0$ we find $T_2>0$ such that for any $\epsi\in(0,1)$, $\kappa\in[-1,1]$ and all $t>T_2$ 
\begin{align*}
\|\uepkap(\cdot,t)\|_{\Lo[2]}^2<\frac{\delta}{2}\quad\text{and}\quad\int_t^{t+1}\big\|\nepkap-\overline{n_0}\big\|_{\Lo[2]}^2<\frac{\delta}{2C_1\|\nabla\phi\|_{\Lo[\infty]}^2}.
\end{align*}
Therefore, for $t>T_2$ integrating \eqref{eq:uepkap-decay-eq1} with respect to time shows 
\begin{align*}
\int_{t}^{t+1}\!\intomega\big|\nabla\uepkap\big|^2\leq \intomega\big|\uepkap(\cdot,t)\big|^2+C_1\|\nabla\phi\|_{\Lo[\infty]}^2\int_t^{t+1}\!\intomega\big(\nepkap-\overline{n_0}\big)^2<\frac{\delta}{2}+\frac{\delta}{2}=\delta,
\end{align*}
proving \eqref{eq:uepkap-nab-decay}. Finally, \eqref{eq:uepkap-decay} is an immediate consequence of \eqref{eq:uepkap-nab-decay} and $W_{0,\sigma}^{1,2}(\Omega)\hookrightarrow\Lo[6]$.
\end{bew}

Making use of semigroup estimates for the Stokes equation we can further refine the smallness results of the previous Lemmas.

\begin{lemma}
\label{lem:improved-uepkap-decay}
For all $\delta>0$ and any $p>3$ there exists $T>0$ such that for each $\epsi\in(0,1)$, $\kappa\in[-1,1]$ and all $t>T$ the solution $(\nepkap,\cepkap,\uepkap)$ of \eqref{approxprop} satisfies
\begin{align*}
\big\|\uepkap(\cdot,s)\big\|_{\Lo[p]}<\delta\quad\text{for any }s\in[t,t+1].
\end{align*}
\end{lemma}

\begin{bew}
This is a consequence of Lemmas \ref{lem:nepkap-decay}, \ref{lem:uepkap-decay} and a fixed point argument relying on the regularizing effects of the Stokes semigroup. The proof we give here is based on \cite[Lemma 7.5]{win_chemonavstokesfinal_TransAm17} and \cite[Lemma 3.8]{Lan-Longterm_M3AS16}. We fix some $p_0\in(3,p)$ satisfying $p_0\leq 6$ and the let $\gamma:=\frac{3}{2}(\frac{1}{p_0}-\frac{1}{p})$, noting that by these choices $\gamma$ fulfills $\gamma\in(0,\frac{1}{2}-\frac{3}{2p})$, and hence the constant
\begin{align*}
C_1:=\int_0^1(1-\sigma)^{-\frac{1}{2}-\frac{3}{2p}}\sigma^{-2\gamma}\intd \sigma
\end{align*}
is finite. Moreover, according to the well known smoothing properties of the Stokes operator (\cite{gig86}) there exist $C_2,C_3,C_4>0$ such that
\begin{align}\label{eq:improved-uepkap-decay1}
\|e^{-tA}\mathcal{P}\varphi\|_{\Lo[p]}&\leq C_2t^{-\gamma}\|\varphi\|_{\Lo[p_0]}&&\text{for all }\varphi\in\Lo[p_0]\text{ and all }t>0,\nonumber\\
\|e^{-tA}\mathcal{P}\varphi\|_{\Lo[p]}&\leq C_3\|\varphi\|_{\Lo[p]}&&\text{for all }\varphi\in\Lo[p]\text{ and all }t>0,\\
\|e^{-tA}\mathcal{P}\nabla\cdot\varphi\|_{\Lo[p]}&\leq C_4 t^{-\frac{1}{2}-\frac{3}{2p}}\|\varphi\|_{\Lo[\frac{p}{2}]}&&\text{for all }\varphi\in\Lo[\frac{p}{2}]\text{ and all }t>0.\nonumber
\end{align}
Now, given $\delta>0$ we next fix $\delta_0\in(0,\delta)$ such that
\begin{align*}
 \delta_0 3^{\frac{1}{2}-\frac{3}{2p}-\gamma}C_1C_4<\frac{1}{3}
\end{align*}
and then, in light of Lemma \ref{lem:uepkap-decay}, Lemma \ref{lem:nepkap-decay} and \eqref{phidef}, pick $T_0>0$ such that for any $\epsi\in(0,1)$, $\kappa\in[-1,1]$ and all $t>T_0$
\begin{align*}
\int_t^{t+1}\!\big\|\uepkap\big\|_{\Lo[p_0]}<\frac{\delta_0}{3C_2}\quad\text{and}\quad\int_{t}^{t+3}\big\|\nepkap-\overline{n_0}\big\|_{\Lo[p]}<\frac{\delta_0}{3^{1+\gamma} C_3\|\nabla\phi\|_{\Lo[\infty]}},
\end{align*}
which in particular also entails that for any fixed $t_1>T_0$ and each $\epsi\in(0,1)$ and $\kappa\in[-1,1]$ there exists $t_\star\in(t_1,t_1+1)$ such that $\big\|\uepkap(\cdot,t_\star)\big\|_{\Lo[p_0]}<\frac{\delta_0}{3C_2}$ holds. Letting $X:=\{\varphi:\Omega\times(t_\star,t_\star+3)\to\R\,\vert\,\|\varphi\|_X:=\sup_{s\in(t_\star,t_\star+3)}(t-t_\star)^{\gamma}\|\varphi(s)\|_{\Lo[p]}<\infty\}$ we now consider the map $\Psi$ acting on the closed subset $S:=\{\varphi\in X\,\vert\,\|\varphi\|_{X}\leq\delta_0\}$ defined by
\begin{align*}
\Psi(\varphi)(\cdot,t):=e^{-(t-t_\star)A}\uepkap(\cdot,t_\star)+\int_{t_\star}^t e^{(t-s)A}\mathcal{P}\left(-\nabla\cdot(Y_\epsi\varphi\otimes\varphi)(\cdot,s)+\nepkap(\cdot,s)\nabla\phi\right)\intd s.
\end{align*}
Drawing on \eqref{eq:improved-uepkap-decay1}, the contraction property of $Y_\epsi$ and the Cauchy-Schwarz inequality,
 we find that
\begin{align*}
\|\Psi(\varphi)(\cdot,t)\|_{\Lo[p]}&\leq C_2(t-t_\star)^{-\gamma}\big\|\uepkap(\cdot,t_\star)\big\|_{\Lo[p_0]}\\
&\qquad+C_4\int_{t_\star}^t(t-s)^{-\frac{1}{2}-\frac{3}{2p}}\|\varphi\|_{\Lo[p]}^2\intd s+C_3\|\nabla\phi\|_{\Lo[\infty]}\int_{t_\star}^{t}\big\|\nepkap(\cdot,t)-\overline{n_0}\big\|_{\Lo[p]}\intd s
\end{align*}
for all $t\in(t_\star,t_\star+3)$. In light of our choice for $\delta_0$, the definition of $S$ and the fact that $|t-t_\star|\leq 3$, this implies that
\begin{align}\label{eq:improved-uepkap-decay2}
(t-t_\star)^{\gamma}\|\Psi(\varphi)(\cdot,t)\|_{\Lo[p]}&\leq C_2\big\|\uepkap(\cdot,t_\star)\big\|_{\Lo[p_0]}+ \delta_0^2(t-t_\star)^{\gamma} C_4\int_{t_\star}^{t}(t-s)^{-\frac{1}{2}-\frac{3}{2p}}s^{-2\gamma}\intd s\nonumber\\&\qquad+3^\gamma C_3\|\nabla\phi\|_{\Lo[\infty]}\int_{t_\star}^{t_\star+3}\big\|\nepkap(\cdot,t)-\overline{n_0}\big\|_{\Lo[p]}\intd s\nonumber\\
&< \frac{\delta_0}{3}+\delta_0^2(t-t_\star)^{\gamma-\frac{1}{2}-\frac{3}{2p}-2\gamma+1} C_4\int_0^1(1-\sigma)^{-\frac{1}{2}-\frac{3}{2p}}\sigma^{-2\gamma}\intd \sigma+\frac{\delta_0}{3}\nonumber\\
&\leq\frac{\delta_0}{3}+\delta_0\left(\delta_0 3^{\frac{1}{2}-\frac{3}{2p}-\gamma} C_1C_4\right)+\frac{\delta_0}{3}<\delta_0,
\end{align} 
and hence $\Psi$ maps $S$ onto itself. Similarly, taking into account that for any $\varphi,\psi\in \Lo[p]$
\begin{align*}
\|\varphi\otimes\varphi-\psi\otimes\psi\|_{\Lo[\frac{p}{2}]}\leq (\|\varphi\|_{\Lo[p]}+\|\psi\|_{\Lo[p]})\|\varphi-\psi\|_{\Lo[p]}
\end{align*}
we find that for any $\varphi,\psi\in S$
\begin{align*}
\|\Psi(\varphi)-\Psi(\psi)\|_{\Lo[p]}&\leq C_4\int_{t_\star}^t(t-s)^{-\frac{1}{2}-\frac{3}{2p}}\|\varphi\otimes\varphi-\psi\otimes\psi\|_{\Lo[\frac{p}{2}]}\intd s\\&\leq 2\delta_0 C_4\int_{t_\star}^t(t-s)^{-\frac{1}{2}-\frac{3}{2p}}s^{-2\gamma}\|\varphi-\psi\|_{X}\quad\text{on }(t_\star,t_\star+3),
\end{align*}
so that 
\begin{align*}
(t-t_\star)^{\gamma}\|\Psi(\varphi)-\Psi(\psi)(\cdot,t)\|_{\Lo[p]}\leq 2\delta_03^{\frac{1}{2}-\frac{3}{2p}-\gamma}C_1C_4\|\varphi-\psi\|_X\quad\text{for all }t\in(t_\star,t_\star+3).	
\end{align*}
Since $2\delta_03^{\frac{1}{2}-\frac{3}{2p}-\gamma}C_1C_4<\frac{2}{3}$, $\Psi:S\to S$ is a contracting map and therefore, there exists a unique fixed point of $\Psi$ on $S$, which has to coincide with $\uepkap$ on $(t_\star,t_\star+3)$ (\cite[Theorem V.2.5.1]{sohr}) and we conclude from \eqref{eq:improved-uepkap-decay2} and the fact that $(t_1+2,t_1+3)\subset(t_\star+1,t_\star+3)$ that for any $\epsi\in(0,1)$, $\kappa\in[-1,1]$
\begin{align*}
\big\|\uepkap(\cdot,t)\big\|_{\Lo[p]}<\delta\quad\text{for all }t\in(t_1+2,t_1+3),
\end{align*}
completing the proof.
\end{bew}

\setcounter{equation}{0} 
\section{Uniform eventual smoothness estimates}\label{sec7:evsmooth}
In order to obtain an improvement on the regularity of our solution components we will incorporate
arguments shown in \cite[Lemmas 3.9, 3.10 and 3.11]{Lan-Longterm_M3AS16}. For this to work we will require the following cut-off
functions (cf. \cite{win_chemonavstokesfinal_TransAm17} and \cite{Lan-Longterm_M3AS16}).

\begin{definition}\label{def:xidef}
Given any monotonically increasing function $\xi_0\in\CSp{\infty}{\R}$ satisfying 
\begin{align*}
0\leq\xi_0\leq 1\ \text{on}\ \R,\quad\xi_0\equiv0\ \text{on}\ (-\infty,0]\quad\text{and}\quad\xi_0\equiv1\ \text{on}\ (1,\infty)
\end{align*}
and some $t_0>0$ we set $$\xi_{t_0}(t):=\xi_0(t-t_0),\quad t\in\R.$$
\end{definition}

Relying on well known maximal Sobolev estimates for the Stokes equation we can a uniform bound for $\uepkap$ in certain Hölder spaces.

\begin{lemma}
\label{lem:hoelder-reg-uepkap}
There exist $\gamma\in(0,1)$, $T>0$ and $C>0$ such that for any $\epsi\in(0,1)$, $\kappa\in[-1,1]$ and all $t>T$ the solution $(\nepkap,\cepkap,\uepkap)$ of \eqref{approxprop} satisfies
\begin{align}\label{eq:hoelder-reg-uepkap0}
\big\|\uepkap\big\|_{\CSp{1+\gamma,\frac{\gamma}{2}}{\bomega\times[t,t+1]}}\leq C.
\end{align}
\end{lemma}

\begin{bew} The proof follows the approach undertaken in \cite[Lemma 3.9]{Lan-Longterm_M3AS16}, which relies on maximal Sobolev regularity properties of the Stokes equation and the uniform bounds already prepared. 

Let us first fix the following parameters. Let $s>3$, $r>1$ and then we pick $s_1>2s$ and $s_1'$ such that $\frac{1}{s_1}+\frac{1}{s_1'}=1$. Then according to Lemma \ref{lem:n-lp-bounds} we can find $T'>0$ and $C_1>0$ such that for any $\epsi\in(0,1)$, $\kappa\in[-1,1]$ and all $t>T'$
\begin{align}\label{eq:hoelder-reg-uepkap1}
\int_t^{t+1}\!\big\|\nepkap\big\|_{\Lo[s]}\leq C_1
\end{align}
holds. Moreover, drawing on Lemma \ref{lem:uepkap-decay} we can fix $T>T'$ such that for any $\epsi\in(0,1)$, $\kappa\in[-1,1]$ and all $t>T$ we also have
\begin{align}\label{eq:hoelder-reg-uepkap2}
\big\|\uepkap\big\|_{\LSp{\infty}{(t,t+2);\Lo[r]}}\leq C_1,\quad \big\|\uepkap\big\|_{\LSp{\infty}{(t,t+2);\Lo[s]}}\leq C_1,\quad \big\|\uepkap\big\|_{\LSp{\infty}{(t,t+2);\Lo[s_1']}}\leq C_1.
\end{align}
Now, for $t_0>T$ we let $\xi:=\xi_{t_0}$ denote the cut-off function given by Definition \ref{def:xidef} and find that $\xi\uepkap$ fulfills
\begin{align*}
\begin{array}{r@{\,}l@{\quad}l}
\big(\xi\uepkap\big)_t&=\Delta\big(\xi\uepkap\big)-\kappa\big(Y_\epsi\uepkap\cdot\nabla\big)\xi\uepkap+\nabla\big(\xi\uepkap\big)+\xi\nepkap\nabla\phi+\xi'\uepkap\ \;\text{in }\Omega\times(t_0,\infty),\\
\nabla\cdot(\xi\uepkap)&=0\ \;\text{in }\Omega\times(t_0,\infty),\\
\text{with}\ \; \big(\xi\uepkap\big)(\cdot,t_0)&=0\ \; \text{in }\Omega,\qquad\text{and}\quad\quad\big(\xi\uepkap\big)=0\ \text{on }\romega\times(t_0,\infty).
\end{array}
\end{align*}
Thus, the maximal Sobolev regularity estimate for the Stokes semigroup (\cite{GigSohr91}) provides $C_2>0$ such that for all $\epsi\in(0,1)$ and $\kappa\in[-1,1]$
\begin{align}\label{eq:hoelder-reg-uepkap3}
&\int_{t_0}^{t_0+2}\!\big\|\big(\xi\uepkap\big)_t\big\|_{\Lo[s]}^s+\int_{t_0}^{t_0+2}\!\big\|D^2\big(\xi\uepkap\big)\big\|_{\Lo[s]}^s\\\leq\ & C_2\Big(0+\int_{t_0}^{t_0+2}\!\big\|\mathcal{P}\big(\big(\xi Y_\epsi\uepkap\cdot\nabla\big)\uepkap)\big\|_{\Lo[s]}^s+\int_{t_0}^{t_0+2}\!\big\|\mathcal{P}\big(\xi\nepkap\nabla\phi\big)\big\|_{\Lo[s]}^s+\int_{t_0}^{t_0+2}\!\big\|\mathcal{P}\big(\xi'\uepkap\big)\big\|_{\Lo[s]}^s\Big)\nonumber.
\end{align}
According to \eqref{eq:hoelder-reg-uepkap1} and \eqref{eq:hoelder-reg-uepkap2} we obtain $C_3>0$ such that for any $\epsi\in(0,1)$ and $\kappa\in[-1,1]$ we may estimate 
\begin{align}\label{eq:hoelder-reg-uepkap4}
\int_{t_0}^{t_0+2}\!\big\|\mathcal{P}\big(\xi\nepkap\nabla\phi\big)\big\|_{\Lo[s]}^s+\int_{t_0}^{t_0+2}\!\big\|\mathcal{P}\big(\xi'\uepkap\big)\big\|_{\Lo[s]}^s\leq C_3.
\end{align}
Moreover, there is $C_4>0$ such that for any $\epsi\in(0,1)$, $\kappa\in[-1,1]$ and all $t_0>T$
\begin{align*}
\big\|\mathcal{P}\big(Y_\epsi\uepkap\cdot\nabla\big)\xi\uepkap\big\|_{\Lo[s]}^s&\leq C_4\big\|Y_\epsi\uepkap\big\|_{\Lo[s_1']}^s\big\|\nabla\big(\xi\uepkap\big)\big\|_{\Lo[s_1]}^s\leq C_4 C_1^s\big\|\nabla\big(\xi\uepkap\big)\big\|_{\Lo[s_1]}^s
\end{align*}
on $(t_0,t_0+2)$, due to Hölder's inequality and the fact that $\|Y_\epsi\varphi\|_{\Lo[s_1']}\leq \|\varphi\|_{\Lo[s_1]}$ holds for all $\varphi\in\Lo[s_1']$. Employing the \GNI\ we then obtain $C_5>0$ such that for any $\epsi\in(0,1)$, $\kappa\in[-1,1]$ and all $t_0>T$
\begin{align*}
\big\|\mathcal{P}\big(Y_\epsi\uepkap\cdot\nabla\big)\xi\uepkap\big\|_{\Lo[s]}^s&\leq C_5 C_4 C_1^s\big\|D^2\big(\xi\uepkap\big)\big\|_{\Lo[s]}^{as}\big\|\xi\uepkap\big\|_{\Lo[r]}^{(1-a)s}\leq C_5 C_4 C_1^{(2-a)s}\big\|D^2\big(\xi\uepkap\big)\big\|_{\Lo[s]}^{as}
\end{align*}
holds on $(t_0,t_0+2)$, where $a=\frac{\frac{1}{3}-\frac{1}{s_1}+\frac{1}{r}}{\frac{2}{3}-\frac{1}{s}+\frac{1}{r}}\in(\frac{1}{2},1)$. Hence, upon integration with respect to time an application of Young's inequality provides $C_6>0$ such that for any $\epsi\in(0,1)$, $\kappa\in[-1,1]$ and all $t_0>T$
\begin{align*}
C_2\int_{t_0}^{t_0+2}\!\big\|\mathcal{P}\big(Y_\epsi\uepkap\cdot\nabla\big)\xi\uepkap\big\|_{\Lo[s]}^s&\leq \frac{1}{2}\int_{t_0}^{t_0+2}\!\big\|D^2\big(\xi\uepkap\big)\big\|_{\Lo[s]}+2 C_2 C_6,
\end{align*}
which combined with \eqref{eq:hoelder-reg-uepkap3} and \eqref{eq:hoelder-reg-uepkap4} shows that for any $\epsi\in(0,1)$, $\kappa\in[-1,1]$ and all $t_0>T$ we have
\begin{align*}
\int_{t_0}^{t_0+2}\!\big\|\big(\xi\uepkap\big)_t\big\|_{\Lo[s]}^s+\frac{1}{2}\int_{t_0}^{t_0+2}\!\big\|D^2\big(\xi\uepkap\big)\big\|_{\Lo[s]}^s\leq 2 C_6C_2+C_3C_2.
\end{align*}
Due to $\xi\equiv 1$ on $(t_0+1,t_0+2)$, this readily implies that for any $s>1$ there exist $C_7>0$ and $T>0$ such that for any $\epsi\in(0,1)$, $\kappa\in[-1,1]$ and all $t>T$ 
\begin{align*}
\int_{t}^{t+1}\!\big\|u_{\epsi t}^{(\kappa)}\big\|_{\Lo[s]}^s+\int_t^{t+1}\!\big\|\uepkap\big\|^s_{\W[2,s]}\leq C_7,
\end{align*}
and in light of known embedding results (e.g. \cite[Theorem 1.1]{Amann00}) entails \eqref{eq:hoelder-reg-uepkap0}.
\end{bew}

Arguments along the same lines of the previous lemma (and previously also employed in \cite[Lemmas 3.10 and 3.11]{Lan-Longterm_M3AS16}), this time drawing on maximal Sobolev estimates for the Neumann heat semigroup, also help us derive Hölder bounds for the remaining components. We proceed with proving a corresponding bound for the signal chemical.

\begin{lemma}
\label{lem:hoelder-reg-cepkap}
For any $p\in(1,\infty)$ there exist $T>0$ and $C>0$ such that for each $\epsi\in(0,1)$, $\kappa\in[-1,1]$ and all $t>T$ the solution $(\nepkap,\cepkap,\uepkap)$ of \eqref{approxprop} satisfies
\begin{align*}
\int_t^{t+1}\big\|c_{\epsi t}^{(\kappa)}\big\|_{\Lo[p]}+\int_t^{t+1}\big\|\cepkap\big\|_{\W[2,p]}\leq C.
\end{align*}
Furthermore, there exist $\gamma\in(0,1)$, $T>0$ and $C'>0$ such that for each $\epsi\in(0,1)$, $\kappa\in[-1,1]$ and all $t>T$
\begin{align}\label{eq:hoelder-reg-cepkap1}
\big\|\cepkap\big\|_{\CSp{1+\gamma,\frac{\gamma}{2}}{\bomega\times[t,t+1]}}\leq C'.
\end{align}
\end{lemma}

\begin{bew}
Given an arbitrary $p\in(1,\infty)$ we first fix $q\in(1,p)$. Now, according to Lemmas \ref{lem:hoelder-reg-uepkap} and \ref{lem:n-lp-bounds} we can pick $T'>0$ and $C_1>0$ such that for all $\epsi\in(0,1)$ and $\kappa\in[-1,1]$ we have
\begin{align*}
\big\|\uepkap\big\|_{\LSp{\infty}{\Omega\times(T',\infty)}}+\big\|\nepkap\big\|_{\LSp{\infty}{(T',\infty);\Lo[p]}}\leq C_1.
\end{align*}
Then, for any $t_0>T'$ we denote by $\xi:=\xi_{t_0}$ a temporal cutoff function as given by Defintion \ref{def:xidef} and observe that $\xi\cepkap$ then satisfies
\begin{align*}
\big(\xi\cepkap\big)_t+\xi\uepkap\cdot\nabla\cepkap=\Delta\big(\xi\cepkap\big)+\frac{1}{\epsi}\xi\cepkap\ln\big(1+\epsi\nepkap\big)+\xi'\cepkap\quad\text{on }\Omega\times(t_0,\infty)
\end{align*}
with $\frac{\partial(\xi\cepkap)}{\partial\nu}=0$ on $\romega\times(t_0,\infty)$ and $\xi\cepkap(\cdot,t_0)=0$ in $\Omega$. In light of the maximal Sobolev regularity estimates for the Neumann heat semigroup (\cite{GigSohr91}), \eqref{eq:approxest} and \eqref{eq:locex-bounds} this implies the existence of $C_2>0$ such that for all $\epsi\in(0,1)$ and $\kappa\in[-1,1]$
\begin{align}\label{eq:hoelder-reg-cepkap2}
\int_{t_0}^{t_0+2}\!\big\|\big(&\xi\cepkap\big)_t\big\|_{\Lo[p]}^p+\int_{t_0}^{t_0+2}\!\big\|\Delta\big(\xi\cepkap\big)\big\|_{\Lo[p]}^p\nonumber\\
&\leq C_2\Big(0+\int_{t_0}^{t_0+2}\!\big\|\xi\uepkap\cdot\nabla\cepkap\big\|_{\Lo[p]}^p+\int_{t_0}^{t_0+2}\!\big\|\xi\cepkap\nepkap\big\|_{\Lo[p]}^p+\int_{t_0}^{t_0+2}\!\big\|\xi'\cepkap\big\|_{\Lo[p]}^p\Big)\\
&\leq C_2C_1^p\int_{t_0}^{t_0+2}\big\|\nabla\big(\xi\cepkap\big)\big\|_{\Lo[p]}^p+2C_2C_1^p\|c_0\|_{\Lo[\infty]}^p+2C_2\|c_0\|_{\Lo[\infty]}^p\|\xi'\|_{\LSp{\infty}{\R}}^p\nonumber
\end{align}
holds. As the \GNI\ entails the existence of $C_3>0$ such that for all $\epsi\in(0,1)$, $\kappa\in[-1,1]$ and $t>t_0$ we have
\begin{align*}
\big\|\nabla\big(\xi\cepkap\big)(\cdot,t)\big\|_{\Lo[p]}^p&\leq C_3\big\|\Delta\big(\xi\cepkap\big)(\cdot,t)\big\|_{\Lo[p]}^{ap}\big\|\xi\cepkap(\cdot,t)\big\|_{\Lo[q]}^{(1-a)p}+C_3\big\|\xi\cepkap(\cdot,t)\big\|_{\Lo[\infty]}^p,
\end{align*}
with $a=\frac{\frac{1}{3}-\frac{1}{p}+\frac{1}{q}}{\frac{2}{3}-\frac{1}{p}+\frac{1}{q}}$ satisfying $a\in(\frac{1}{2},1)$, an employment of Young's inequality together with Hölder's inequality and \eqref{eq:locex-bounds} provides $C_4>0$ such that for all $\epsi\in(0,1)$, $\kappa\in[-1,1]$ and any $t_0>T'$ the inequality from \eqref{eq:hoelder-reg-cepkap2} reads like
\begin{align*}
\int_{t_0}^{t_0+2}\!\big\|\big(\xi\cepkap\big)\big\|_{\Lo[p]}^p&+\frac{1}{2}\int_{t_0}^{t_0+2}\!\big\|\Delta\big(\xi\cepkap\big)\big\|_{\Lo[p]}^p\\&\leq \|c_0\|_{\Lo[\infty]}^p(2C_4+2C_3C_2C_1^p+2C_2C_1^p+2C_2\|\xi'\|_{\LSp{\infty}{\R}}^p).
\end{align*}
Since $\xi\equiv 1$ on $(t_0+1,t_0+2)$, this shows that for any $p>1$ one can find $C_5>0$ and $T:=T'+1>0$ such that for any $t>T$ and all $\epsi\in(0,1)$, $\kappa\in[-1,1]$ we have
\begin{align*}
\int_{t}^{t+1}\!\big\|c_{\epsi t}^{(\kappa)}\big\|_{\Lo[p]}^p+\int_{t}^{t+1}\!\big\|\cepkap\big\|_{\W[2,p]}\leq C_5.
\end{align*}
The asserted Hölder regularity finally results from an application of an embedding result e.g. presented in \cite[Theorem 1.1]{Amann00} by taking $p$ large enough.
\end{bew}

A final iteration of similar arguments entails a uniform Hölder bound for the first solution component.

\begin{lemma}
\label{lem:hoelder-reg-nepkap}
There exist $\gamma\in(0,1)$, $T>0$ and $C>0$ such that for each $\epsi\in(0,1)$, $\kappa\in[-1,1]$ and all $t>T$ the solution
$(\nepkap,\cepkap,\uepkap)$ of \eqref{approxprop} satisfies
\begin{align}\label{eq:hoelder-reg-nepkap1}
\big\|\nepkap\big\|_{\CSp{1+\gamma,\frac{\gamma}{2}}{\bomega\times[t,t+1]}}\leq C.
\end{align}
\end{lemma}

\begin{bew}
We work along similar lines as in the previous lemma. First, given any $p>1$ we pick $q\in(1,p)$ and, in light of Lemmas \ref{lem:n-lp-bounds}, \ref{lem:hoelder-reg-uepkap} and \ref{lem:hoelder-reg-cepkap}, can then find $T'>0$ and $C_1>0$ such that for any $\epsi\in(0,1)$, $\kappa\in[-1,1]$ we have
\begin{align*}
\big\|\nepkap\big\|_{\LSp{\infty}{(T',\infty);\Lo[p]}}+\big\|\nepkap\big\|_{\LSp{\infty}{(T',\infty);\Lo[q]}}+\big\|\nepkap\big\|_{\LSp{\infty}{(T',\infty);\Lo[2p]}}\leq C_1,
\end{align*}
and $C_2,C_3>0$ such that for all $t>T'$ the estimates
\begin{align*}
&\big\|\nabla\cepkap(\cdot,t)\big\|_{\LSp{\infty}{\Omega\times(t,t+2)}}\leq C_2,\qquad\int_t^{t+2}\big\|\Delta\cepkap(\cdot,t)\big\|_{\Lo[2p]}^{p}\leq C_2,\quad\text{and}\quad\big\|\uepkap\big\|_{\LSp{\infty}{\Omega\times(t,t+2)}}\leq C_3
\end{align*}
hold. Now, for $t_0>T'$ we once more denote by $\xi:=\xi_{t_0}$ the cutoff function from Definition \ref{def:xidef} and the maximal Sobolev regularity estimates (\cite{GigSohr91}) then again entail the existence of $C_4>0$ such that for all $\epsi\in(0,1)$ and $\kappa\in[-1,1]$
\begin{align}\label{eq:hoelder-reg-nepkap-eq1}
&\int_{t_0}^{t_0+2}\!\big\|\big(\xi n^{(\kappa)}_{\epsi}\big)_t\big\|_{\Lo[p]}^p+\int_{t_0}^{t_0+2}\!\big\|\Delta\big(\xi\nepkap\big)\big\|_{\Lo[p]}^p\\&\ \leq C_4\int_{t_0}^{t_0+2}\!\big\|\xi\uepkap\cdot\nabla\nepkap\big\|_{\Lo[p]}^p+C_4\int_{t_0}^{t_0+2}\bigg\|\xi\nabla\cdot\bigg(\frac{\nepkap}{1+\epsi\nepkap}\nabla\cepkap\bigg)\bigg\|_{\Lo[p]}^p+C_4\int_{t_0}^{t_0+2}\!\big\|\xi'\nepkap\big\|_{\Lo[p]}^p.\nonumber
\end{align}
Next, to estimate mixed derivative term we note that by the bounds prepared at the start of the lemma
\begin{align}\label{eq:hoelder-reg-nepkap-eq2}
\bigg\|\xi\nabla\cdot\bigg(\frac{\nepkap}{1+\epsi\nepkap}\nabla\cepkap\bigg)\bigg\|_{\Lo[p]}^p&=\bigg\|\xi\frac{\nabla\nepkap\cdot\nabla\cepkap}{(1+\epsi\nepkap)^2}+\xi\frac{\nepkap}{1+\epsi\nepkap}\Delta\cepkap\bigg\|_{\Lo[p]}^p\nonumber\\
&\leq 2^p C_2^p\big\|\nabla\big(\xi\nepkap\big)\big\|_{\Lo[p]}^p+2^p C_1^p\big\|\Delta\cepkap\big\|_{\Lo[2p]}^p
\end{align}
is valid on $(t_0,t_0+2)$. Moreover, the \GNI\ implies the existence of $C_5>0$ such that
\begin{align*}
\|\nabla\varphi\|_{\Lo[p]}^p\leq C_5\|\Delta\varphi\|_{\Lo[p]}^{ap}\|\varphi\|_{\Lo[q]}^{(1-a)p}+\|\varphi\|_{\Lo[q]}^p\quad\text{for all }\varphi\in\W[2,p],
\end{align*}
where $a=\frac{\frac{1}{3}+\frac{1}{q}-\frac{1}{p}}{\frac{2}{3}+\frac{1}{q}-\frac{1}{p}}\in(\frac12,1)$, and hence we infer from Young's inequality that there is $C_6>0$ such that
\begin{align}\label{eq:hoelder-reg-nepkap-eq3}
C_4(C_3^p+2^p C_2^p)\|\nabla\varphi\|_{\Lo[p]}^p\leq \frac12\|\Delta\varphi\|_{\Lo[p]}^p+C_6\|\varphi\|_{\Lo[q]}^p\quad\text{for all }\varphi\in\W[2,p].
\end{align}
Thus, collecting \eqref{eq:hoelder-reg-nepkap-eq1}--\eqref{eq:hoelder-reg-nepkap-eq3}, we conclude for all $\epsi\in(0,1)$ and $\kappa\in[-1,1]$
\begin{align*}
\int_{t_0}^{t_0+2}&\!\big\|\big(\xi n^{(\kappa)}_{\epsi}\big)_t\big\|_{\Lo[p]}^p+\int_{t_0}^{t_0+2}\!\big\|\Delta\big(\xi\nepkap\big)\big\|_{\Lo[p]}^p\\&\ \leq C_4(C_3^p+2^p C_2^p)\int_{t_0}^{t_0+2}\!\big\|\nabla\big(\xi\nepkap\big)\big\|_{\Lo[p]}^p+2^p C_1^p\int_{t_0}^{t_0+2}\!\big\|\Delta\cepkap\big\|_{\Lo[2p]}^p+2C_1^pC_4\|\xi'\|_{\LSp{\infty}{\R}}\\
&\ \leq \frac{1}{2}\int_{t_0}^{t_0+2}\!\big\|\Delta\big(\xi\nepkap\big)\big\|_{\Lo[p]}^p+2C_6 C_1^p+2^{p} C_1^p C_2+2C_1^pC_4\|\xi'\|_{\LSp{\infty}{\R}},
\end{align*}
which, due to $\xi\equiv 1$ on $(t_0+1,t_0+2)$ implies the existence of $C_7>0$ such that for any $\epsi\in(0,1)$, $\kappa\in[-1,1]$ and all $t>T:=T'+1$ we have
\begin{align*}
\int_{t}^{t+1}&\!\big\| n^{(\kappa)}_{\epsi t}\big\|_{\Lo[p]}^p+\frac{1}{2}\int_{t}^{t+1}\!\big\|\Delta\nepkap\big\|_{\Lo[p]}^p\leq C_7.
\end{align*}
Taking $p$ large enough, the desired Hölder regularity is again an immediate consequence of the embedding result in e.g. \cite[Theorem 1.1]{Amann00}.
\end{bew}

In light of standard parabolic theory and the Arzelà-Ascoli theorem  we can make use of the uniform estimates from the previous three lemmas to conclude that after an eventual smoothing time $\Td>0$ the solution obtained in Lemma \ref{lem:epsi-limit} is actually a classical solution.

\begin{lemma}
\label{lem:ev-smooth}
There exist $\gamma>0$ and $\Td>0$ such that for each $\kappa\in[-1,1]$ the weak solution $(\nkap,\ckap,\ukap)$ of \eqref{CTNScons} obtained in Lemma \ref{lem:epsi-limit} satisfies
\begin{align*}
\nkap,\ckap\in\CSp{2+\gamma,1+\frac{\gamma}{2}}{\bomega\times[\Td,\infty)}\quad\text{and}\quad\ukap\in\CSp{2+\gamma,1+\frac{\gamma}{2}}{\bomega\times[\Td,\infty);\R^3}.
\end{align*}
In particular, $(\nkap,\ckap,\ukap)$ together with some $P^{(\kappa)}\in\CSp{1,0}{\bomega\times(\Td,\infty)}$ solve \eqref{CTNScons} classically in $\Omega\times(\Td,\infty)$. Moreover, there exists $C>0$ such that for all $\kappa\in[-1,1]$ and all $t\geq\Td$
\begin{align}\label{eq:kappa-indep-C2-bound}
\big\|\nkap(\cdot,t)\big\|_{\CSpnl{2+\gamma,1+\frac{\gamma}{2}}{\bomega\times[t,t+1]}}+\big\|\ckap(\cdot,t)\big\|_{\CSpnl{2+\gamma,1+\frac{\gamma}{2}}{\bomega\times[t,t+1]}}+\big\|\ukap(\cdot,t)\big\|_{\CSpnl{2+\gamma,1+\frac{\gamma}{2}}{\bomega\times[t,t+1]}}\leq C.
\end{align}
\end{lemma}

\begin{bew}
We employ standard parabolic regularity theory in a similar fashion as e.g. displayed in \cite[Lemma 3.12]{Lan-Longterm_M3AS16}. Drawing on Lemmas \ref{lem:hoelder-reg-uepkap}, \ref{lem:hoelder-reg-cepkap} and \ref{lem:hoelder-reg-nepkap} we can pick some $\gamma'\in(0,1)$, $T>0$ and $C>0$ such that \eqref{eq:hoelder-reg-uepkap1}, \eqref{eq:hoelder-reg-cepkap1} and \eqref{eq:hoelder-reg-nepkap1} hold for any $t>T$, each $\epsi\in(0,1)$ and every $\kappa\in[-1,1]$. Accordingly, by making use of the Arzelà-Ascoli theorem we have
\begin{align*}
\nepkap\to\nkap,\ \cepkap\to\ckap\quad\text{in}\ \;\CSp{1+\gamma',\frac{\gamma'}{2}}{\bomega\times[t,t+1]}\ \;\text{and}\ \;\uepkap\to\ukap\quad\text{in}\ \;\CSp{1+\gamma',\frac{\gamma'}{2}}{\bomega\times[t,t+1];\R^3}
\end{align*}
along a subsequence of the sequence $(\epsi_j)_{j\in\N}$ obtained in Lemma \ref{lem:epsi-limit}, the members of which, for convenience, we still label $\epsi_j$. Now, letting $\xi:=\xi_T$ be given by Definition \ref{def:xidef} we note that $\xi\ckap$ solves
\begin{align*}
c_t=\Delta c + g,\quad c(T)=0,\quad\frac{\partial c}{\partial\nu}\big|_{\romega}=0,
\end{align*}
in the weak sense with $g=-\xi\nkap\ckap-\xi\ukap\nabla\ckap+\ckap\xi'\in\CSp{\gamma'}{\bomega\times(T,\infty)}$. In light of standard parabolic theory (e.g. \cite[III.5.1 and IV.5.3]{lsu}) we can hence conclude that for some $\gamma_1\in(0,1)$ $\ckap\in\CSp{2+\gamma_1,1+\frac{\gamma_1}{2}}{\bomega\times[T+1,\infty)}$ and that for $\gamma\leq\gamma_1$ there is $C_1>0$ such that for any $\kappa\in[-1,1]$ \eqref{eq:kappa-indep-C2-bound} is true for $\ckap$. In a similar fashion we observe that $\xi\nkap$ is a solution of
\begin{align*}
n_t=\Delta n- a\cdot\nabla n+ b,\quad n(T)=0,\quad\frac{\partial n}{\partial\nu}\big|_{\romega}=0,
\end{align*}
with $a=\nabla\ckap+\ukap$ and $b=-\xi\nkap\Delta c+\nkap\xi'$ both being of class $\CSp{\gamma',\frac{\gamma'}{2}}{\bomega\times(T,\infty)}$ and employing parabolic regularity theory (e.g. \cite[IV.5.3 and III.5.1]{lsu}) once more, we find $\gamma_2\in(0,1)$ such that $\nkap\in\CSp{2+\gamma_2,1+\frac{\gamma_2}{2}}{\bomega\times[T+1,\infty)}$ and that for $\gamma\leq\gamma_2$ there is $C_2>0$ such that \eqref{eq:kappa-indep-C2-bound} is valid for $\nkap$. Lastly, since $\xi\ukap$ solves
\begin{align*}
u_t=\Delta u+h:=\mathcal{P}\big(\xi' \ukap-\kappa\xi(\ukap\cdot\nabla)\ukap)+\xi\nkap\nabla\phi\big),\quad\nabla\cdot u=0,\quad u(T-1)=0,\quad u\big|_{\romega}=0,
\end{align*}
where $h$ is again Hölder continuous due to the bounds from Lemmas \ref{lem:hoelder-reg-uepkap}, \ref{lem:hoelder-reg-cepkap} and \ref{lem:hoelder-reg-nepkap} and \eqref{phidef}. 
Hence, the Schauder theory for Stokes equation (e.g. \cite[Theorem 1.1]{Solonnikov-SchauderEst07}) combined with the uniqueness property (\cite[V.1.5.1]{sohr}) entails that for some $\gamma_3\in(0,1)$ we have $\ukap\in\CSp{2+\gamma_3,1+\frac{\gamma_3}{2}}{\bomega\times[T+1,\infty)}$ and that for $\gamma\leq\gamma_3$ \eqref{eq:kappa-indep-C2-bound} is also valid for $\ukap$. Letting $\gamma:=\min\{\gamma_1,\gamma_2,\gamma_3\}$ we obtain the inclusion in the asserted function spaces, whereas the existence of a corresponding $P^{(\kappa)}\in\CSp{1,0}{\bomega\times(\Td,\infty)}$ such that $(\nkap,\ckap,\ukap,P^{(\kappa)})$ solves \eqref{CTNScons} classically in $\Omega\times(\Td,\infty)$ is an immediate consequence of these regularity properties (\cite{sohr}).
\end{bew}

\setcounter{equation}{0} 
\section{Uniform exponential decay after the smoothing time}\label{sec8:expdec}
For the remainder of the work we will denote by $\Td>0$ the smoothing time obtained in Lemma \ref{lem:ev-smooth}.
Employing a second Aubing-Lions type argument for taking $\kappa\to0$, we are still left with the obstacle that this limit procedure will only yield convergence on compact subsets of $\bomega\times[0,\infty)$. In order to extend the convergence beyond compact subsets our next objective will be to improve the previously obtained stabilization properties to a more detailed decay including an exponential rate of convergence, which, on the one hand, will still be independent of $\kappa\in[-1,1]$ and, on the other, will be valid for all $t>\Td$. We start with supplementing our decay results by the following lemma.  

\begin{lemma}
\label{lem:dec-n-infty}
For all $\delta>0$ one can find $T\geq \Td$ such that for each $\kappa\in[-1,1]$ and all $t>T$ the solution $(\nkap,\ckap,\ukap)$ of \eqref{CTNScons} satisfies
\begin{align*}
\big\|\nkap(\cdot,t)-\overline{n_0}\big\|_{\Lo[\infty]}<\delta.
\end{align*}
\end{lemma}

\begin{bew}
According to the \GNI\ there is $C_1>0$ such that for any $\kappa\in[-1,1]$ and all $t>\Td$
\begin{align}\label{eq:dec-n-infty-eq1}
\big\|\nkap(\cdot,t)-\overline{n_0}\big\|_{\Lo[\infty]}^5&\leq C_1\big\|\nkap(\cdot,t)-\overline{n_0}\big\|_{\W[1,\infty]}^3 \big\|\nkap(\cdot,t)-\overline{n_0}\big\|_{\Lo[2]}^2.
\end{align}
Moreover, by \eqref{eq:kappa-indep-C2-bound} and \eqref{IR} we can find $C_2>0$ such that for any $\kappa\in[-1,1]$ and $t>\Td$
\begin{align}\label{eq:dec-n-infty-eq2}
C_1\big\|\nkap(\cdot,t)-\overline{n_0}\big\|_{\W[1,\infty]}^3 \leq C_1\big(\big\|\nkap(\cdot,t)\big\|_{\CSp{1}{\bomega}}+\|\overline{n_0}\|_{\Lo[\infty]}\big)^3\leq C_2,
\end{align}
since $\overline{n_0}$ is spatially homogeneous. Then, given $\delta>0$ we set $\delta_0:=\frac{\delta^5}{C_2}$ and rely on Lemma \ref{lem:nepkap-decay} to find $T>\Td$ such that for any $\kappa\in[-1,1]$ and $t>T$
\begin{align*}
\big\|\nkap(\cdot,t)-\overline{n_0}\big\|_{\Lo[2]}^2<\delta_0.
\end{align*}
A combination of this with \eqref{eq:dec-n-infty-eq1} and \eqref{eq:dec-n-infty-eq2} yields that for any $\kappa\in[-1,1]$ and all $t>T$
\begin{align*}
\big\|\nkap(\cdot,t)-\overline{n_0}\big\|_{\Lo[\infty]}<(C_2\delta_0)^\frac{1}{5}=\delta
\end{align*}
holds, finalizing the proof.
\end{bew}

Combining the previous lemma with the fact that $(\nkap,\ckap,\ukap)$ solves \eqref{CTNScons} classically on $\Omega\times(\Td,\infty)$, we can improve the eventual decay of the oxygen, which in Lemma \ref{lem:decay-c-Linfty} was still of a quite general nature, to a decay with exponential rate.

\begin{lemma}
\label{lem:exp-dec-c}
There exist $\mu>0$ and $C>0$ such that for each $\kappa\in[-1,1]$ and all $t>0$ the solution $(\nkap,\ckap,\ukap)$ of \eqref{CTNScons} satisfies
\begin{align}\label{eq:exp-dec-c-linfty}
\big\|\ckap(\cdot,t)\big\|_{\Lo[\infty]}\leq Ce^{-\mu t}.
\end{align}
Moreover, for any $p\geq 1$ there exist $\mu'>0$ and $C'>0$ such that for each $\kappa\in[-1,1]$ and all $t>\Td$
\begin{align}\label{eq:exp-dec-c-W1p}
\big\|\ckap(\cdot,t)\big\|_{\W[1,p]}\leq C'e^{-\mu' t}.
\end{align}
\end{lemma}

\begin{bew} We follow the reasoning of \cite[Lemmas 4.5 and 4.6]{WangWinklerXiang-smallconvection-MathZ18}. Drawing on the $\kappa$--independent stabilization property obtained in Lemma \ref{lem:dec-n-infty} we can fix $T>\Td$ such that
\begin{align*}
\nkap\geq C_1:=\frac{\overline{n_0}}{2}\quad\text{in }\Omega\times(T,\infty),
\end{align*}
where $C_1$ is positive due to \eqref{IR}. Noting that $(\nkap,\ckap,\ukap)$ is a classical solution of \eqref{CTNScons} on $\Omega\times(\Td,\infty)$, we can make use of the second equation of \eqref{CTNScons} to find that
\begin{align*}
c^{(\kappa)}_t\leq \Delta\ckap-\ukap\cdot\nabla\ckap-C_1\ckap\quad\text{in }\Omega\times(T,\infty),
\end{align*}
and therefore, the comparison principle combined with \eqref{eq:locex-bounds} implies that
\begin{align*}
\ckap(\cdot,t)\leq\big\|\ckap(\cdot,T)\big\|_{\Lo[\infty]}e^{-C_1(t-T)}\leq\|c_0\|_{\Lo[\infty]}e^{-C_1(t-T)}\quad\text{for all }t>T.
\end{align*}
Relying once more on \eqref{eq:locex-bounds}, we find that \eqref{eq:exp-dec-c-linfty} also holds for $0<t\leq T$ by letting $C:=\|c_0\|_{\Lo[\infty]}e^{C_1 T}$ and $\mu:=\frac{\overline{n_0}}{2}$.
As for the decay involving the gradient, we note that, assuming $p>3$, the \GNI\ provides $C_2>0$ such that
\begin{align*}
\big\|\ckap(\cdot,t)\big\|_{\W[1,p]}\leq C_2\big\|\ckap(\cdot,t)\big\|_{\CSpnl{2}{\bomega}}^\frac{p-3}{2p}\big\|\ckap(\cdot,t)\big\|_{\Lo[\infty]}^\frac{p+3}{2p}
\end{align*}
is valid for all $t>\Td$, which according to \eqref{eq:exp-dec-c-linfty} and \eqref{eq:kappa-indep-C2-bound} implies \eqref{eq:exp-dec-c-W1p}. 
\end{bew}

With the previous result at hand, we can not only transfer the exponential rate of convergence to the first solution component, but also establish this decay starting from the smoothing time $\Td$, clarifying the convergence statement from Lemma \ref{lem:dec-n-infty}.

\begin{lemma}
\label{lem:exp-dec-n}
There exist $\mu>0$ and $C>0$ such that for each $\kappa\in[-1,1]$ and all $t>\Td$ the solution $(\nkap,\ckap,\ukap)$ of \eqref{CTNScons} satisfies
\begin{align}\label{eq:exp-dec-n-linfty}
\big\|\nkap(\cdot,t)-\overline{n_0}\big\|_{\Lo[\infty]}< Ce^{-\mu t}.
\end{align}
Moreover, for any $p\geq 1$ there exist $\mu'>0$ and $C'>0$ such that for each $\kappa\in[-1,1]$ and all $t>\Td$
\begin{align*}
\big\|\nkap(\cdot,t)-\overline{n_0}\big\|_{\W[1,p]}\leq C'e^{-\mu' t}.
\end{align*}
\end{lemma}

\begin{bew}
We adjust the arguments of \cite[Lemma 4.7]{WangWinklerXiang-smallconvection-MathZ18} to our setting and start by working along similar lines as in Lemma \ref{lem:nepkap-decay}, while this time making sure we keep the $L^2$ norm of $\nabla \cepkap$ to make full use of the exponential decay established in \ref{lem:exp-dec-c}. In fact, drawing on the first equation in \eqref{CTNScons}, as well as integration by parts, Young's inequality and the Poincaré inequality we obtain $C_1>0$ such that for any $\kappa\in[-1,1]$ and all $t>\Td$
\begin{align*}
\frac{\intd}{\intd t}\intomega\big(\nkap-\overline{n_0}\big)^2
\leq -\frac{1}{C_1}\intomega\big(\nkap-\overline{n_0}\big)^2+\sup_{\kappa'\in[-1,1]}\big\|n^{(\kappa')}\big\|_{\LSp{\infty}{\Omega\times(\Td,\infty)}}\intomega\big|\nabla\ckap\big|^2
\end{align*}
holds. Hence, according to \eqref{eq:kappa-indep-C2-bound} and Lemma \ref{lem:exp-dec-c}, we can fix $\mu_1>0$ with $\frac{1}{C_1}>\mu_1$ such that for any $\kappa\in[-1,1]$ the function $y^{(\kappa)}(t):=\intomega(\nkap(\cdot,t)-\overline{n_0})^2$ satisfies
$\frac{\intd}{\intd t}y^{(\kappa)}(t)+\frac{1}{C_1}y^{(\kappa)}(t)\leq C_2e^{-\mu_1 t}$ for all $t>\Td$, which implies
\begin{align*}
y^{(\kappa)}(t)\leq \left(C_3+\frac{C_1C_2}{1-C_1\mu_1}\right)e^{\mu_1\Td}e^{-\mu_1 t}\quad\text{for all }t>\Td,
\end{align*}
with $C_3:=\intomega\big(\nkap(\cdot,\Td)-\overline{n_0}\big)^2$ being independent of $\kappa$ and finite, again due to \eqref{eq:kappa-indep-C2-bound} and \eqref{IR}. Interpolation using the \GNI\ and, once more, \eqref{eq:kappa-indep-C2-bound} finally extends to \eqref{eq:exp-dec-n-linfty} upon appropriately adjusting the constants. For the decay of the Sobolev-Norm, we assume, again without loss of generality, that $p>3$ and draw on the \GNI\ to find $C_4>0$ such that
\begin{align*}
\big\|\nkap(\cdot,t)-\overline{n_0}\big\|_{\W[1,p]}&\leq C_4\big\|\nkap(\cdot,t)-\overline{n_0}\big\|_{\W[2,\infty]}^\frac{p-3}{2p}\big\|\nkap(\cdot,t)-\overline{n_0}\big\|_{\Lo[\infty]}^\frac{p+3}{2p}\\
&\leq C_4\Big(\big\|\nkap(\cdot,t)\big\|_{\CSpnl{2}{\bomega}}+\big\|\overline{n_0}\big\|_{\Lo[\infty]}\Big)^\frac{p-3}{2p}\big\|\nkap(\cdot,t)-\overline{n_0}\big\|_{\Lo[\infty]}^\frac{p+3}{2p}\quad\text{for all }t>\Td,
\end{align*}
because $\overline{n_0}$ is constant in space. Hence, the claimed exponential decay is a consequence of \eqref{eq:exp-dec-n-linfty}, \eqref{IR} and Lemma \ref{lem:ev-smooth}.
\end{bew}

In the final part of this section, we extend the exponential stabilization of the first component to the fluid velocity field.	

\begin{lemma}
\label{lem:exp-dec-u}
There exist $\mu>0$ and $C>0$ such that for each $\kappa\in[-1,1]$ and all $t>\Td$ the solution $(\nkap,\ckap,\ukap)$ of \eqref{CTNScons} satisfies
\begin{align}\label{eq:exp-dec-u-linfty}
\big\|\ukap(\cdot,t)\big\|_{\Lo[\infty]}< Ce^{-\mu t}.
\end{align}
Moreover, for any $p\geq 1$ there exists $C'>0$ such that for each $\kappa\in[-1,1]$ and all $t>\Td$
\begin{align*}
\big\|\ukap(\cdot,t)\big\|_{\W[1,p]}\leq C'e^{-\mu t}.
\end{align*}
\end{lemma}

\begin{bew}
Similar to the previous two lemmas and inspired by \cite[Lemma 4.8]{WangWinklerXiang-smallconvection-MathZ18}, we proceed to derive an exponential decay estimate for the fluid velocity. Due to $\nabla\cdot\ukap=0$ in $\Omega\times(0,\infty)$ and $\ukap=0$ on $\romega\times(0,\infty)$ we obtain upon testing the third equation of \eqref{CTNScons} against $\ukap$ that for any $\kappa\in[-1,1]$ and all $t>\Td$ 
\begin{align}\label{eq:exp-dec-u-eq1}
\frac{1}{2}\frac{\intd}{\intd t}\intomega\big|\ukap\big|^2+\intomega\big|\nabla\ukap\big|^2=\intomega\big(\nkap-\overline{n_0}\big)\nabla\phi\cdot\ukap.
\end{align}
Since Poincaré's inequality provides $C_1>0$ such that for any $\kappa\in[-1,1]$ and all $t>\Td$ we have
\begin{align*}
C_1 \intomega\big|\ukap\big|^2\leq \intomega\big|\nabla\ukap\big|^2,
\end{align*}
we can employ the Hölder and Young inequalities to conclude from \eqref{eq:exp-dec-u-eq1} that for any $\kappa\in[-1,1]$ and all $t>\Td$
\begin{align*}
\frac{1}{2}\frac{\intd}{\intd t}\intomega\big|\ukap\big|^2+
\frac{C_1}{2}\intomega\big|\ukap\big|^2+\frac{1}{2}\intomega\big|\nabla\ukap\big|^2&\leq \sqrt{|\Omega|}\|\nabla\phi\|_{\Lo[\infty]}\big\|\nkap-\overline{n_0}\big\|_{\Lo[\infty]}\big\|\ukap\big\|_{\Lo[2]}\\
&\leq \frac{\sqrt{|\Omega|}}{\sqrt{C_1}}\|\nabla\phi\|_{\Lo[\infty]}\big\|\nkap-\overline{n_0}\big\|_{\Lo[\infty]}\big\|\nabla\ukap\big\|_{\Lo[2]}\\
&\leq C_2\big\|\nkap-\overline{n_0}\big\|_{\Lo[\infty]}^2+\frac{1}{2}\intomega\big|\nabla\ukap\big|^2,
\end{align*}
where $C_2:=\frac{|\Omega|\|\nabla\phi\|_{\Lo[\infty]}^2}{2C_1}$. Hence, making use the decay estimate from Lemma \ref{lem:exp-dec-n}, we can find $\mu_1\in(0,C_1)$ and $C_3>0$ such that $y^{(\kappa)}(t):=\intomega\big|\ukap(\cdot,t)\big|^2$, $t>\Td$ satisfies
\begin{align*}
\frac{\intd}{\intd t}y^{(\kappa)}(t)+C_1 y^{(\kappa)}(t)\leq C_3e^{-\mu_1 t}\quad\text{for all }t>\Td,
\end{align*}
implying that for any $\kappa\in[-1,1]$ and all $t>\Td$
\begin{align}\label{eq:exp-dec-u-eq2}
y^{(\kappa)}(t)\leq \left(y^{(\kappa)}(\Td)+\frac{C_3}{C_1-\mu}\right)e^{C_1 \Td}e^{-\mu_1 t}=: C_4 e^{-\mu_1 t},
\end{align}
where $C_4=(y(\Td)+\frac{C_3}{C_1-\mu})e^{C_1 \Td}$ does not depend on $\kappa$ and is finite due to \eqref{eq:kappa-indep-C2-bound}. Now, with $\alpha\in(\frac{3}{4},1)$ given by \eqref{IR} and according to \cite[Theorem 1.4.4]{hen81} there is $C_5>0$ such that for any $\kappa\in[-1,1]$ and each $t>\Td$
\begin{align*}
\big\|A^\alpha\ukap(\cdot,t)\big\|_{\Lo[2]}\leq C_5\big\|A\ukap(\cdot,t)\big\|_{\Lo[2]}^\alpha\big\|\ukap(\cdot,t)\big\|_{\Lo[2]}^{1-\alpha}\leq C_5\big\|\ukap(\cdot,t)\big\|_{\CSp{2}{\bomega}}^\alpha\big\|\ukap(\cdot,t)\big\|_{\Lo[2]}^{1-\alpha},
\end{align*}
and drawing once more on \eqref{eq:kappa-indep-C2-bound} and \eqref{eq:exp-dec-u-eq2} together with the embedding $D(A^\alpha)\hookrightarrow\LSp{\infty}{\Omega;\R^3}$ (\cite[Theorem 1.6.1]{hen81}) provides $C_6>0$ such that for any $\kappa\in[-1,1]$ and all $t>\Td$
\begin{align*}
\big\|\ukap(\cdot,t)\big\|_{\Lo[\infty]}\leq C_6 e^{-\mu_1(1-\alpha)t}
\end{align*}
holds and hence proves \eqref{eq:exp-dec-u-linfty}. Employing the \GNI\ in a similar fashion as in the proofs of the previous two lemmas finally entails the exponential decay of the desired Sobolev norms.
\end{bew}

\setcounter{equation}{0} 
\section{The second limit. Taking \texorpdfstring{$\kappa\to 0$}{kappa to zero}}\label{sec9:kappalimit}
The uniform exponential decay starting from the smoothing time $\Td$ was the last missing ingredient for proving our theorem. Before we give the proof of the theorem however, we first collect many of the prepared estimates for the following second limit procedure.
\begin{lemma}
\label{lem:kappa-limit}
Given any null sequence $(\kappa_j)_{j\in\N}\subset[-1,1]$ one can find a subsequence $(\kappa_{j_k})_{k\in\N}$ and functions 
\begin{align*}
n&\in\LSploc{\frac53}{\bomega\times[0,\infty)}\quad\text{with}\quad\nabla n\in\LSploc{\frac54}{\bomega\times[0,\infty)},\\
c&\in\LSp{\infty}{\Omega\times(0,\infty)}\quad\text{with}\quad\nabla c\in\LSploc{4}{\bomega\times[0,\infty)},\\
u&\in L^2_{loc}\big([0,\infty);W_{0,\sigma}^{1,2}(\Omega)\big),
\end{align*}
such that the global weak solution $(\nkap,\ckap,\ukap)$ of \eqref{CTNScons} satisfies
\begin{alignat*}{2}
\nkap&\to n \qquad&&\text{in }\LSploc{p}{\bomega\times[0,\infty)}\text{ for any }p\in[1,\tfrac53)\text{ and a.e. in }\Omega\times(0,\infty),\\
\nabla \nkap&\wto \nabla n &&\text{in }\LSploc{\frac54}{\bomega\times[0,\infty)},\\
\nkap&\wto n&&\text{in }\LSploc{\frac53}{\bomega\times[0,\infty)},\\
\ckap&\to c  &&\text{in }\LSploc{p}{\bomega\times[0,\infty)}\text{ for any }p\in[1,\infty)\text{ and a.e. in }\Omega\times(0,\infty),\\
\ckap&\wsto c&&\text{in }\LSp{\infty}{\Omega\times(0,\infty)},\\
\nabla\ckap &\wto\nabla c &&\text{in }\LSploc{4}{\bomega\times[0,\infty)},\\
\ukap&\to u &&\text{in }\LSploc{2}{\bomega\times[0,\infty)}\text{ and a.e. in }\Omega\times(0,\infty),\\
\ukap&\wto u &&\text{in }\LSploc{\frac{10}3}{\bomega\times[0,\infty)},\\
\nabla\ukap&\wto\nabla u &&\text{in }\LSploc{2}{\bomega\times[0,\infty)},
\end{alignat*}
as $\kappa=\kappa_{j_k}\to0$. The triple $(n,c,u)$ is a global weak solution of the chemotaxis-Stokes system \eqref{CTStokes}, \eqref{BC} and \eqref{IC} in the sense of Definition \ref{def:sol}, and one can find $P\in\CSp{1,0}{\bomega\times(\Td,\infty)}$ such that $(n,c,u,P)$ solve \eqref{CTStokes}, \eqref{BC}, \eqref{IC} classically in $\Omega\times(\Td,\infty)$. Moreover, there exist $\mu>0$ and $C>0$ such that for all $t>\Td$,
\begin{align}\label{eq:exp-dec-stokes1}
\|n(\cdot,t)-\overline{n_0}\|_{\Lo[\infty]}+\|c(\cdot,t)\|_{\Lo[\infty]}+\|u(\cdot,t)\|_{\Lo[\infty]}<C e^{-\mu t}
\end{align}
and for any $p\geq1$ there are $\mu'>0$ and $C'>0$ such that
\begin{align}\label{eq:exp-dec-stokes2}
\|n(\cdot,t)-\overline{n_0}\|_{\W[1,p]}+\|c(\cdot,t)\|_{\W[1,p]}+\|u(\cdot,t)\|_{\W[1,p]}<C' e^{-\mu' t}
\end{align}
is valid for all $t>\Td$.
\end{lemma}

\begin{bew}
As the bounds in Lemmas \ref{lem:bounds}, \ref{lem:timespace-bounds} and \ref{lem:timereg} are independent of $\epsi\in(0,1)$ and $\kappa\in[-1,1]$ they are inherited by the limit functions $\nkap,\ckap$ and $\ukap$ obtained in Lemma \ref{lem:epsi-limit} and hence an identical reasoning, drawing on the Aubin--Lions Lemma \cite[Corollary 8.4]{Sim87} and Vitali's theorem, as previously done in Lemma \ref{lem:epsi-limit}, establishes the asserted convergence properties and weak solution properties of the limit functions $n,c$ and $u$. That there exists some $P\in\CSp{1,0}{\bomega\times(\Td,\infty)}$, which together with $(n,c,u)$ solves \eqref{CTStokes} classically in $\Omega\times(\Td,\infty)$ is then a consequence of Lemma \ref{lem:ev-smooth} and the Arzelà--Ascoli theorem. The exponential decay estimates for times larger than the smoothing time $\Td$, as stated in \eqref{eq:exp-dec-stokes1} and \eqref{eq:exp-dec-stokes2}, are a consequence of Lemmas \ref{lem:exp-dec-c}, \ref{lem:exp-dec-n} and \ref{lem:exp-dec-u}.
\end{bew}

With the limit objects and local convergence properties prepared by the previous lemma, can finally draw on the uniform exponential decay for large times established in Section \ref{sec8:expdec} to extend the local convergence to convergence beyond compact subsets of $\bomega\times[0,\infty)$, as claimed in the main theorem.

\begin{proof}[\textbf{Proof of Theorem \ref{theo:1}}:]
The existence and the regularity and solution properties of the claimed functions were already established in Lemmas \ref{lem:epsi-limit} and \ref{lem:kappa-limit}. We are left with verifying the convergence with respect to the desired norms as in \eqref{eq:thm1}. According to Lemma \ref{lem:exp-dec-n} and \eqref{eq:exp-dec-stokes1}, given any $p_1\in[1,\frac{5}{3})$ we can fix $\mu>0$ and $C_1>0$ such that for any $\kappa\in[-1,1]$
\begin{align}\label{eq:theo-proof0}
\big\|n^{(\kappa)}(\cdot,t)&-\overline{n_0}\big\|_{\Lo[p_1]}^{p_1}+\|n(\cdot,t)-\overline{n_0}\|_{\Lo[p_1]}^{p_1}<C_1 e^{-\mu t}
\end{align}
holds for all $t>\Td$. Now, given $\delta>0$ we pick $T_\star\geq\max\big\{\Td,\frac1{p_1\mu}\ln\big(\frac{2C_1}{p_1\mu\delta}\big)\big\}$ and obtain from \eqref{eq:theo-proof0} that
\begin{align}\label{eq:theo-proof1}
\big\|n^{(\kappa)}-n\big\|_{\LSp{p_1}{(T_\star,\infty),\Lo[p_1]}}^{p_1}&\leq\int_{T_\star}^\infty\!\big\|n^{(\kappa)}(\cdot,t)-\overline{n_0}\big\|_{\Lo[p_1]}^{p_1}\intd t+\int_{T_\star}^\infty\!\big\|n(\cdot,t)-\overline{n_0}\big\|_{\Lo[p_1]}^{p_1}\intd t\nonumber\\&\leq 2C_1\int_{T_\star}^\infty e^{-p_1\mu t}\intd t=\frac{2C_1}{p_1\mu}e^{-p_1 \mu T_\star}\leq \frac{\delta}{2}\quad\text{for any }\kappa\in[-1,1].
\end{align} 
Next, given a null sequence $(\kappa_j)_{j\in\N}\subset[-1,1]$ and denoting by $(\kappa_{j_k})_{k\in\N}$ the subsequence from Lemma \ref{lem:kappa-limit}, we can conclude from the convergence statements in Lemma \ref{lem:kappa-limit} that
\begin{align*}
n^{(\kappa_{j_k})}\to n\quad\text{in }\LSploc{p_1}{[0,\infty),\Lo[p_1]}\quad\text{as }\kappa_{j_k}\to0.
\end{align*}
Hence, for the given $\delta>0$ there is some $k_0\in\N$ such that
\begin{align}\label{eq:theo-proof2}
\big\|n^{(\kappa_{j_k})}-n\big\|_{\LSp{p_1}{[0,T_\star];\Lo[p_1]}}^{p_1}\leq\frac{\delta}{2}\quad\text{is valid for all }k\geq k_0. 
\end{align}
Combination of \eqref{eq:theo-proof1} and \eqref{eq:theo-proof2} shows that for all $k\geq k_0$ we have
\begin{align*}
\big\|n^{(\kappa_{j_k})}-n\big\|_{\LSp{p_1}{[0,\infty);\Lo[p_1]}}^{p_1}=\big\|n^{(\kappa_{j_k})}-n\big\|_{\LSp{p_1}{[0,T_\star];\Lo[p_1]}}^{p_1}+\big\|n^{(\kappa_{j_k})}-n\big\|_{\LSp{p_1}{(T_\star,\infty),\Lo[p_1]}}^{p_1}\leq \delta,
\end{align*}
from which we conclude the first part of \eqref{eq:thm1}. Similarly, drawing on Lemma \ref{lem:exp-dec-n} and \eqref{eq:exp-dec-stokes2}, for given $p_2\in[1,\frac{5}{4})$ we can fix $\mu'>0$ and $C_2>0$ such that for any $\kappa\in[-1,1]$
\begin{align*}
\big\|\nkap(\cdot,t)-\overline{n_0}\big\|_{\W[1,p_2]}^{p_2}+\|n(\cdot,t)-\overline{n_0}\|_{\W[1,p_2]}^{p_2}\leq C_2 e^{-\mu't}
\end{align*}
holds for all $t>\Td$, from which we once again conclude that for the given $\delta>0$ we can pick $T_\star'\geq\max\{\Td,\frac{1}{p_2\mu'}\ln(\frac{2C_2}{p_2\mu'\delta})\}$ such that
\begin{align}\label{eq:theo-proof3}
\big\|\nabla\nkap-\nabla n\big\|_{\LSp{p_2}{(T_\star',\infty);\Lo[p_2]}}^{p_2}\leq\frac{\delta}{2}\quad\text{for any }\kappa\in[-1,1].
\end{align}
Since we know from Lemma \ref{lem:kappa-limit} that
\begin{align*}
\nabla n^{(\kappa_{j_k})}\wto\nabla n\quad\text{in }\LSploc{p_2}{[0,\infty);\Lo[p_2]}\quad\text{as }\kappa_{j_k}\to0,
\end{align*}
we can make use of the fact that $p_2<\frac{5}{4}$ to employ Vitali's theorem in combination with the uniform bounds presented in Lemma \ref{lem:timespace-bounds} to find that actually
\begin{align*}
\nabla n^{(\kappa_{j_k})}\to\nabla n\quad\text{in }\LSploc{p_2}{[0,\infty);\Lo[p_2]}\quad\text{as }\kappa_{j_k}\to0.
\end{align*}
From this we conclude that there is some $k_0'\in\N$ such that
\begin{align}\label{eq:theo-proof4}
\big\|\nabla n^{(\kappa_{j_k})}-\nabla n\big\|_{\LSp{p_2}{[0,T_\star'];\Lo[p_2]}}^{p_2}\leq\frac{\delta}{2}\quad\text{holds for all }k\geq k_0',
\end{align}
so that a combination of \eqref{eq:theo-proof3} and \eqref{eq:theo-proof4} again entails the convergence in the desired topology. Analogous arguments drawing on Lemma \ref{lem:exp-dec-c}, Lemma \ref{lem:exp-dec-u}, \eqref{eq:exp-dec-stokes1}, \eqref{eq:exp-dec-stokes2} and the uniform bounds in the Lemmas \ref{lem:bounds} and \ref{lem:timespace-bounds} finally entail the remaining properties listed in \eqref{eq:thm1}, completing the proof. 
\end{proof}

\section*{Acknowledgements}
The author acknowledges support of the {\em Deutsche Forschungsgemeinschaft} in the context of the project
  {\em Emergence of structures and advantages in cross-diffusion systems (project no. 411007140)}.

\footnotesize{
\setlength{\bibsep}{2pt plus 0.5ex}

}
\end{document}